\documentclass[12pt]{article}

\usepackage{indentfirst}
\usepackage{amsfonts}
\usepackage{amsmath}
\usepackage{amssymb}
\usepackage{amsbsy, amsthm}
\usepackage[margin=1in]{geometry}

\newtheorem{dfn}{Definition}[section]
\newtheorem{theorem}[dfn]{Theorem}
\newtheorem{lemma}[dfn]{Lemma}
\newtheorem{proposition}[dfn]{Proposition}

\newtheorem{corollary}[dfn]{Corollary}
\newtheorem{question}[dfn]{Question}
\newenvironment{pf}{\noindent{\bf Proof.}}
{\enspace\vrule height5pt depth0pt width5pt}

\def\T {{\mathcal T}}

\def\F {{\mathcal F}}

\def\D {{\mathcal D}}

\def\C {{\mathcal C}}

\def\V {{\mathcal V}}
\def\LL{{\mathcal L}}

\begin{document}
\title{Proper conflict-free list-coloring, odd minors, subdivisions, and layered treewidth}
\author{Chun-Hung Liu\thanks{chliu@math.tamu.edu. Partially supported by NSF under award DMS-1954054 and CAREER award DMS-2144042.} \\
\small Department of Mathematics, \\
\small Texas A\&M University,\\
\small College Station, TX 77843-3368, USA}

\maketitle

\begin{abstract}
Proper conflict-free coloring is an intermediate notion between proper coloring of a graph and proper coloring of its square.
It is a proper coloring such that for every non-isolated vertex, there exists a color appearing exactly once in its (open) neighborhood.
Typical examples of graphs with large proper conflict-free chromatic number include graphs with large chromatic number and bipartite graphs isomorphic to the $1$-subdivision of graphs with large chromatic number.
In this paper, we prove two rough converse statements that hold even in the list-coloring setting.
The first is for sparse graphs: for every graph $H$, there exists an integer $c_H$ such that every graph with no subdivision of $H$ is (properly) conflict-free $c_H$-choosable.
The second applies to dense graphs: every graph with large conflict-free choice number either contains a large complete graph as an odd minor or contains a bipartite induced subgraph that has large conflict-free choice number.
These give two incomparable (partial) answers of a question of Caro, Petru\v{s}evski and \v{S}krekovski.
We also prove quantitatively better bounds for minor-closed families, implying some known results about proper conflict-free coloring and odd coloring in the literature.
Moreover, we prove that every graph with layered treewidth at most $w$ is (properly) conflict-free $(8w-1)$-choosable.
This result applies to $(g,k)$-planar graphs, which are graphs whose coloring problems have attracted attention recently.
\end{abstract}

\section{Introduction}

Graph coloring is a central research area in graph theory.
For an integer $k$, a {\it $k$-coloring} of a graph $G$ is a function $\phi: V(G) \rightarrow [k]$.
A coloring of a graph $G$ is {\it proper} if it is a function $\phi$ with domain $V(G)$ such that $\phi(x) \neq \phi(y)$ for every $xy \in E(G)$.
The {\it chromatic number} of $G$, denoted by $\chi(G)$, is the minimum $k$ such that there exists a proper $k$-coloring of $G$.

One topic in graph coloring is about the chromatic number of $G^2$, where $G^2$ is the graph with the same vertex-set as a graph $G$ and two vertices are adjacent in $G^2$ if and only if the distance between them in $G$ is at most 2.
For example, Wegner \cite{w} proposed a conjecture about the chromatic number of $G^2$ for planar graphs $G$ with maximum degree $\Delta$; Erd\H{o}s and Ne\v{s}et\v{r}il (see \cite{en,fgst}) proposed a conjecture about the chromatic number of $(L(G))^2$ for a graph $G$ with maximum degree $\Delta$, where $L(G)$ is the line graph of $G$.
A lot of work about these two conjectures has been done.
For example, see \cite{a_square,bpp_square,bj_square,c_square,hjt_square,h_square,hqt_square,hsy_square,hdk,mr_square,t_square}.
We remark that this is far from a complete list of results about coloring the square of a graph, even just for the aforementioned two conjectures.

Note that $\chi(G^2)$ is tied to the maximum degree of $G$, since the vertices in the neighborhood of any vertex require all different colors.
In this paper we consider proper (open) conflict-free colorings.
Roughly speaking, such a coloring only requires some color appears exactly once in the neighborhood of any vertex, so it is a relaxation of a proper coloring of the square of a graph.
We will define this notion formally later in this paper and observe that the number of required colors is no longer tied to the maximum degree.

Besides being a relaxation of proper colorings of the square of a graph, conflict-free coloring was motivated by a frequency assignment problem for cellular networks \cite{ahkms}.
Even, Lotker, Ron and Smorodinsky \cite{elrs} introduced conflict-free coloring for hypergraphs, which is a vertex-coloring such that for every hyperedge $e$, there exists a color that appears exactly once on the vertices contained in $e$.
Cheilaris \cite{c} considered the special case when the hyperedges are exactly the (open) neighborhoods of the vertices of a graph, while Abel, Alvarez, Demaine, Fekete, Gour, Hesterberg, Keldenich and Scheffer \cite{aadfghks} considered the special case when the hyperedges are exactly the closed neighborhoods of the vertices of a graph.

In this paper, a conflict-free coloring refers to a conflict-free coloring with respect to (open) neighborhoods.
For a vertex $v$ of a graph $G$, we define $N_G(v)=\{u \in V(G): uv \in E(G)\}$ and define $N_G[v]=N_G(v) \cup \{v\}$.

\begin{dfn}
{\rm A {\it conflict-free coloring} of a graph $G$ is a function $\phi$ with domain $V(G)$ such that for every vertex $v$ of $G$, either $N_G(v)=\emptyset$, or there exists an element $c_v$ in the image of $\phi$ such that $|\phi^{-1}(\{c_v\}) \cap N_G(v)|=1$.
That is, if $v$ is not an isolated vertex, then $c_v$ is a color that appears in the neighborhood of $v$ exactly once.}
\end{dfn}

By combining the two notions of proper colorings and (not necessarily proper) colorings and the two notions of conflict-free colorings (with respect to open neighborhoods) and conflict-free colorings with respect to closed neighborhoods, there are four conflict-free-types of colorings studied in the literature (for example, \cite{aadfghks,cps,flrs,pt}).
For a graph $G$, we define $\chi_{{\rm pcf}}(G)$ and $\chi_{{\rm pcfc}}(G)$ to be the minimum $k$ such that $G$ admits a proper conflict-free $k$-coloring and a proper conflict-free $k$-coloring with respect to closed neighborhoods, respectively; we define $\chi_{{\rm icf}}(G)$ and $\chi_{{\rm icfc}}(G)$ to be the minimum $k$ such that $G$ admits a conflict-free $k$-coloring and a conflict-free $k$-coloring with respect to closed neighborhoods, respectively.

It is easy to see that $\chi_{{\rm icf}}(G) \leq \chi_{{\rm pcf}}(G) \geq \chi(G)$ and $\chi_{{\rm icfc}}(G) \leq \chi_{{\rm pcfc}}(G) \geq \chi(G)$.
In fact, for every proper coloring, the color on each vertex $v$ appears exactly once in the closed neighborhood of $v$.
So every proper coloring is a conflict-free coloring with respect to closed neighborhoods.
Hence $\chi_{{\rm pcfc}}(G)=\chi(G)$.
Therefore, $\chi_{{\rm pcf}}(G)$ is an upper bound for the other three parameters $\chi_{{\rm pcfc}}(G)$, $\chi_{{\rm icf}}(G)$ and $\chi_{{\rm icfc}}(G)$.
In addition, as observed by Pach and Tardos \cite{pt}, $\chi_{{\rm icfc}}(G) \leq 2\chi_{{\rm icf}}(G)$.

On the other hand, the gap between $\chi_{{\rm pcf}}(G)$ and the other 3 parameters can be arbitrarily large.
If $G$ is a $1$-subdivision\footnote{The $1$-subdivision of a graph $H$ is the graph obtained from $H$ by subdividing each edge exactly once.} of a graph $H$, then the neighbors of any vertex of degree 2 in $G$ receive different colors in any proper conflict-free coloring, so $\chi_{{\rm pcf}}(G) \geq \chi(H)$, \footnote{An $(\leq k)$-subdivision of a graph $H$ is a graph obtained from $H$ by subdividing each edge $e$ of $H$ $s_e$ times, for some $0 \leq s_e \leq k$. The same observation shows that $\chi_{{\rm pcf}}(G) \geq \chi(H)$ for every graph $G$ that is an $(\leq 1)$-subdivision of $H$.} but $\chi(G)=2$ since $G$ is bipartite.
In addition, for the complete graph $K_n$ on $n \geq 3$ vertices, by coloring one vertex with color 1, one vertex with color 2 and all other vertices color 3, we obtain an improper conflict-free coloring (with respect to both open neighborhoods and closed neighborhoods), so $\chi_{{\rm icf}}(K_n)\leq 3$ and $\chi_{{\rm icfc}}(K_n) \leq 3$, but $\chi_{{\rm pcf}}(K_n) \geq \chi(K_n)=n$.
This implies that having bounded $\chi_{{\rm icf}}$ and $\chi_{{\rm icfc}}$ does not imply having bounded chromatic number.
So even though the analog of Hadwiger's conjecture for conflict-free coloring with respect to closed neighborhoods is true \cite{aadfghks}, it does not give any upper bound for the chromatic number of proper minor-closed families.

In this paper we are interested in finding upper bounds for $\chi_{{\rm pcf}}(G)$. 
The first reason is that proper conflict-free coloring is an intermediate notion between proper coloring of a graph $G$ and proper coloring of its square $G^2$.
Unlike the chromatic number of $G^2$, $\chi_{{\rm pcf}}(G)$ is not necessarily tied to the maximum degree of $G$. 
The second reason is that proper conflict-free coloring behaves very differently from proper coloring. 
As we discussed earlier, 1-subdivision of a graph with large chromatic number has large $\chi_{{\rm pcf}}(G)$.
So having maximum average degree at most 4 does not imply that $\chi_{{\rm pcf}}(G)$ is bounded, even though it ensures bounded chromatic number.
The third reason is that $\chi_{{\rm pcf}}(G)$ is an upper bound of the other three conflict-free-types of parameters $\chi_{{\rm pcfc}}(G)$, $\chi_{{\rm icf}}(G)$ and $\chi_{{\rm icfc}}(G)$.
So it suffices to find upper bounds for $\chi_{{\rm pcf}}(G)$.

\subsection{Our results}

We start with a natural question from \cite{cps}.

\begin{question}[{{\cite[Question 6.2]{cps}}}] \label{que:generic}
Find a ``generic'' graph family ${\mathcal F}$ for which there exists a constant $c$ such that $\chi_{{\rm pcf}}(G)/\chi(G) \leq c$ for every $G \in {\mathcal F}$.
\end{question}

As we mentioned earlier, unlike many coloring parameters, having bounded degeneracy and bounded chromatic number does not ensure an upper bound for proper conflict-free chromatic number.
For every graph $H$, the 1-subdivision of $H$ is 2-degenerate and bipartite, but its proper conflict-free chromatic number is at least the chromatic number of $H$. 
Hence, if a graph class has bounded proper conflict-free chromatic number, then there exist a graph $H_1$ and a bipartite graph $H_2$ such that this class does not contain the $1$-subdivision of $H_1$ and does not contain the graph $H_2$.
That is, 1-subdivisions and bipartiteness are natural obstructions for having small proper conflict-free chromatic number.
Two of our main results of this paper (Theorems \ref{thm:topo_intro} and \ref{thm:odd_minor_intro}) are rough converse of this observation and give partial answers for Question \ref{que:generic}, where one addresses ``subdivisions'' and the other addresses ``induced bipartiteness'', even for the list-coloring setting.

Let $G$ be a graph.
Let $k$ be a real number.
A {\it $k$-list-assignment} of $G$ is a function that maps each vertex of $G$ to a set with size at least $k$.
For a $k$-list-assignment $L$ of $G$, an {\it $L$-coloring} is a function $\phi$ with domain $V(G)$ such that $\phi(v) \in L(v)$ for every $v \in V(G)$.
An $L$-coloring $\phi$ is {\it proper} if $\phi(x) \neq \phi(y)$ for every $xy \in E(G)$.
An $L$-coloring $\phi$ is {\it conflict-free} if for every vertex $v$ of $G$, either $N_G(v)=\emptyset$, or there exists a color $c_v$ in the image of $\phi$ such that $|\phi^{-1}(\{c_v\}) \cap N_G(v)|=1$.
We say that $G$ is {\it conflict-free $k$-choosable}\footnote{We use the terminology ``conflict-free $k$-choosable'' instead of ``proper conflict-free $k$-choosable'' because in the context of graph coloring, being $k$-choosable already requires the corresponding colorings being proper. We follow the same convention here.} if for every $k$-list-assignment $L$ of $G$, there exists a proper conflict-free $L$-coloring of $G$.

Our first answer for Question \ref{que:generic} is the following.

\begin{theorem} \label{thm:topo_intro}
For every graph $H$, there exists a real number $c_H$ such that every graph that does not contain a subdivision of $H$ (as a subgraph) is conflict-free $c_H$-choosable.\footnote{This paper provides two proofs of this result. The first proof is explicitly stated in this paper and relies on the machinery about clique-sums developed in Section \ref{sec:sum}. The second proof is implicitly stated in this paper and uses Lemma \ref{lemma:ordering} proved in Section \ref{sec:layer}. The motivation of Lemma \ref{lemma:ordering} in this paper is to study $(g,k)$-planar graphs and graphs with bounded layered treewidth, which will be described later in this section. The author did not notice the relationship between Lemma \ref{lemma:ordering} and Theorem \ref{thm:topo_intro} until one day before he announced the first draft of this paper on arXiv \cite{l} when a draft of Hickingbotham \cite{h} appeared on arXiv. Hickingbotham \cite{h} independently discovered Lemma \ref{lemma:ordering} with essentially the same proof as ours and observed that known results in the literature immediately show that Lemma \ref{lemma:ordering} applies to graph classes with bounded expansion and hence implies Theorem \ref{thm:topo_intro}. See the concluding remarks in Section \ref{sec:remarks} for details and our further generalization of Hickingbotham's observation. We still keep our original proof of Theorem \ref{thm:topo_intro} by using clique-sums, because this proof is just a simple application of the machinery. This machinery will also be used to prove our second answer for Question \ref{que:generic} in terms of odd minors and induced bipartite subgraphs. In particular, this machinery is applicable to graphs with unbounded degeneracy (so beyond the scope of Lemma \ref{lemma:ordering}) and can probably be further developed to provide better bounds for other notions of colorings.}
\end{theorem}

Note that graphs satisfying Theorem \ref{thm:topo_intro} are sparse in the sense that the number of edges is at most a linear function of the number of vertices.
Our second answer for Question \ref{que:generic} addresses induced bipartite subgraphs and odd minors, which applies to graphs that have quadratically many edges.
We discuss our other results for sparse graphs before stating this second answer. 

Minor-closed families are special cases of Theorem \ref{thm:topo_intro}.
A graph $H$ is a {\it minor} of another graph $G$ if $H$ is isomorphic to a graph that is obtained from a subgraph of $G$ by contracting edges.
A class ${\mathcal F}$ of graphs is {\it minor-closed} if every minor of a member of ${\mathcal F}$ is in ${\mathcal F}$.
A minor-closed family is {\it proper} if it does not contain all graphs.
Clearly, if ${\mathcal F}$ is a proper minor-closed family, then there exists a graph $H$ such that every graph in ${\mathcal F}$ does not contain $H$ as a minor, so every graph in ${\mathcal F}$ does not contain a subdivision of $H$.

We prove a more explicit upper bound for minor-closed families by using degeneracy.
Let $d$ be a real number.
We say a graph $G$ is {\it $d$-degenerate} if every subgraph of $G$ has a vertex of degree at most $d$.
We say a class ${\mathcal F}$ of graphs is {\it $d$-degenerate} if every graph in ${\mathcal F}$ is $d$-degenerate.

A classical result of Mader \cite{m_minor} implies that for every proper minor-closed family ${\mathcal F}$, there exists a real number $k$ such that ${\mathcal F}$ is $k$-degenerate.
The optimal degeneracy for many minor-closed families has been determined exactly or asymptotically, such as for the class of $H$-minor free graphs when $H$ is a small complete graph \cite{j,m_2,st}, a large complete graph \cite{t}, a complete bipartite graph \cite{crs,kp,kp_2,kp_3,ko}, the Petersen graph \cite{hw}, a dense graph \cite{mt,tw}, or a properly 4-colorable graph in a monotone class admitting strongly sublinear balanced separators \cite{hnw}.

We prove that the conflict-free choice number of a proper minor-closed family can be upper bounded in terms of its degeneracy and hence can be combined with the aforementioned results about degeneracy.

\begin{theorem} \label{thm:minor_intro}
Let $d$ be a nonnegative integer.
If $\F$ is a $d$-degenerate minor-closed family, then every graph $G$ in $\F$ is conflict-free $(2d+1)$-choosable.\footnote{We remark that the strength of this result is on the quantitative side. Even though Lemma \ref{lemma:ordering} and the aforementioned independent work of Hickingbotham \cite{h} can be applied to minor-closed families, they provide quantitatively much worse bound than Theorem \ref{thm:minor_intro}, unless some known results in the literature are significantly improved.} 
\end{theorem}

Besides combining Theorem \ref{thm:minor_intro} with aforementioned results about degeneracy, Theorem \ref{thm:minor_intro} already implies a number of results in the literature.

If $\F$ is the class of graphs with treewidth\footnote{Treewidth will be defined in Section \ref{sec:sum}.} at most $w$, then $\F$ is $w$-degenerate, so every graph in $\F$ is conflict-free $(2w+1)$-choosable.
Outerplanar graphs have treewidth at most 2, so Theorem \ref{thm:minor_intro} implies that every outerplanar graph is conflict-free 5-choosable.
It implies an earlier result in \cite{flrs} stating that every outerplanar graph is properly conflict-free 5-colorable, which is tight since the 5-cycle is not properly conflict-free 4-colorable.
In addition, every forest has treewidth at most 1, so Theorem \ref{thm:minor_intro} implies that every forest is conflict-free 3-choosable.
It implies an earlier result in \cite{cps} stating that every forest is properly conflict-free 3-colorable, which is tight since the 3-vertex path is not properly conflict-free 2-colorable.

We remark that it is easy to prove that every graph with maximum degree $\Delta$ is conflict-free $(2\Delta+1)$-choosable. 
See \cite{cl}. 
Moreover, the coefficient 2 can be improved to 1 if $\Delta=3$ \cite{ly} and to 1.6550826 if $\Delta$ is sufficiently large \cite{cl}.
Graphs with maximum degree $\Delta$ are $\Delta$-degenerate, but $1$-degenerate graphs can have arbitrarily large maximum degree.
Theorem \ref{thm:minor_intro} shows that the condition on maximum degree (for coefficient 2) can be replaced by the one on degeneracy if we restrict the graphs to be in minor-closed families.

Another implication of Theorem \ref{thm:minor_intro} is about odd coloring.
A $k$-coloring $\phi$ of a graph $G$ is {\it odd} if for every $v \in V(G)$, either $N_G(v)=\emptyset$, or there exists $c_v$ in the image of $\phi$ such that $|\phi^{-1}(\{c_v\}) \cap N_G(v)|$ is odd.
Clearly, every proper conflict-free coloring is a proper odd coloring.
Cranston, Lafferty and Song \cite{cls} proved that every $d$-degenerate proper minor-closed family is properly odd $(2d+1)$-colorable, and hence this is a special case of Theorem \ref{thm:minor_intro}.

An extensively studied minor-closed family is the class of graphs embeddable in a fixed surface.
We can prove a slightly better bound than Theorem \ref{thm:minor_intro} in this case.

\begin{theorem} \label{thm:surface_intro}
Let $\Sigma$ be a surface with Euler genus $\rho$.
If $\rho \in \{0,1\}$, then every graph embeddable in $\Sigma$ is conflict-free $11$-choosable.
If $\rho \geq 2$, then every graph embeddable in $\Sigma$ is conflict-free $\frac{13+\sqrt{73+48\rho}}{2}$-choosable.
\end{theorem}

Another result of Cranston, Lafferty and Song \cite{cls} is about odd coloring on 1-planar graphs.
For nonnegative integers $g$ and $k$, a graph is {\it $(g,k)$-planar} if it can be drawn in a surface of Euler genus $g$ such that every edge contains at most $k$ crossings.
Note that $(0,k)$-planar graphs are also called {\it $k$-planar graphs} in the literature.
The class of 1-planar graphs is not a topological minor-closed family, so Theorem \ref{thm:topo_intro} does not apply to this class.
On the other hand, Cranston, Lafferty and Song \cite{cls} proved that every 1-planar graph is properly odd 23-colorable.
They \cite{cls} also stated that ``it seems non-trivial to prove a more general result for $k$-planar graphs''.
Another result of this paper solves this case via layered treewidth\footnote{Layered treewidth will be defined in Section \ref{sec:layer}. The author of the present paper (via private communication with Zi-Xia Song) observed that the result in \cite{cls} about graphs with bounded treewidth easily leads to an $O(w)$ upper bound for the proper odd chromatic number of graphs with layered treewidth at most $w$, and hence leads to an $O(k)$ upper bound for the proper odd chromatic number for $k$-planar graphs. Dujmovi\'{c}, Morin and Odak \cite{dmo} later announced a proof for an upper bound $O(k^5)$ for $k$-planar graphs by using strong products of graphs.}.

\begin{theorem} \label{thm:ltw_intro}
If $w$ is a positive integer, then every graph with layered treewidth at most $w$ is conflict-free $(8w-1)$-choosable.
\end{theorem}

It is known \cite{dew} that every $(g,k)$-planar graph has layered treewidth at most $(4g+6)(k+1)$.
Hence Theorem \ref{thm:ltw_intro} gives the following corollary.

\begin{corollary}
For any nonnegative integers $g$ and $k$, every $(g,k)$-planar graph is conflict-free $((32g+48)(k+1)-1)$-choosable.
\end{corollary}

We remark that Hickingbotham \cite{h} independently announced a paper on arXiv when the writing of the first version of this paper \cite{l} was about to be completed. 
The main result in \cite{h} is essentially equivalent to (but actually slightly weaker than) Lemma \ref{lemma:ordering} in this paper that we will develop for proving Theorem \ref{thm:ltw_intro}, with essentially the same proof.
Hickingbotham \cite{h} observed that combining his version of Lemma \ref{lemma:ordering} with a known result in the literature immediately implies that every graph class with bounded expansion has bounded proper conflict-free chromatic number and hence implies Theorem \ref{thm:topo_intro}.
We should address that our proof of Theorem \ref{thm:topo_intro} is a simple application of a machinery about clique-sums, which is part of the main technical contribution of this paper and is used to prove our results about odd minors (Theorem \ref{thm:odd_minor_intro}) that will be described later in this section.
Theorem \ref{thm:odd_minor_intro} applies to dense graphs, so it cannot be proved via Lemma \ref{lemma:ordering} or Hickingbotham's work \cite{h}.
The proof of Theorem \ref{thm:topo_intro} via Lemma \ref{lemma:ordering} or Hickingbotham's work \cite{h} is conceptually simpler than the one via clique-sums, but probably gives weaker quantitative bound.
Moreover, even though Hickingbotham's observation \cite{h} also leads to a $O(w)$ bound for Theorem \ref{thm:ltw_intro}, the coefficient for $w$ in his bound is weaker than the one in our Theorem \ref{thm:ltw_intro}.
Minimizing this coefficient for $w$ is of interest in this paper because a $O(w)$ bound can be easily proved via our Theorem \ref{thm:minor_intro} without resorting to other results in the literature in contrast to Hickingbotham's proof.
As our Lemma \ref{lemma:ordering} is quantitatively stronger than the main result in \cite{h}, our results in this paper with explicit bounds are quantitatively stronger than all results in \cite{h}.
In addition, Theorems \ref{thm:minor_intro} and \ref{thm:surface_intro} are quantitatively better than Lemma \ref{lemma:ordering} and are not covered by Lemma \ref{lemma:ordering} or work in \cite{h}.
In Section \ref{sec:remarks}, we will explain the relationship between our Lemma \ref{lemma:ordering} and Hickingbotham's work \cite{h} in more detail, and we will extend Hickingbotham's observation by combining Lemma \ref{lemma:ordering} with more known results in the literature to immediately give results that are stronger than Hickingbotham's.

Now we address our second partial answer for Question \ref{que:generic}.
Recall that every graph $Q$ that is an $(\leq 1)$-subdivision of a graph $H$ with large chromatic number has large proper conflict-free chromatic number.
It is not hard to show that either $Q$ has large chromatic number, or $Q$ contains an induced subgraph $Q'$ that is a 1-subdivision of a graph with large chromatic number.
(See Proposition \ref{prop:at_most_to_exact} for a proof.)
Note that $Q'$ is an induced bipartite subgraph that has large proper conflict-free chromatic number.
Our second partial answer for Question \ref{que:generic} states that such an induced bipartite graph with large proper conflict-free chromatic number is the only obstruction in odd minor-closed families.

Let $G$ be a graph.
An {\it odd contraction} is the operation that first takes a partition $\{A,B\}$ of $V(G)$ with size two, and then for each connected component $C$ of $G'$, contracts $C$ into a vertex, where $G'$ is the spanning subgraph of $G$ whose edges are exactly the edges of $G$ between $A$ and $B$.
It is easy to show that applying an arbitrary odd contraction on a bipartite graph leads to a bipartite graph.
A graph $H$ is an {\it odd minor} of another graph $G$ if $H$ is isomorphic to a graph that can be obtained from a subgraph of $G$ by repeatedly applying odd contractions.
So every odd minor of a bipartite graph is bipartite.

If $H$ is an odd minor of $G$, then $H$ is a minor of $G$; but the converse is not necessarily true.
The key feature is that unlike $H$-minor free graphs, which have bounded degeneracy, odd $H$-minor free graphs can be dense.
For example, the complete bipartite graph $K_{n,n}$ has $n^2$ edges but does not contain $K_3$ as an odd minor. 

Odd minors were considered when generalizing Hadwiger's conjecture.
Gerards and Seymour (see \cite[p.\ 115]{jt}) conjectured that every graph with no odd $K_{t+1}$-minor is properly $t$-colorable for every positive integer $t$.
It is known \cite{ggrsv} that for every graph $H$, if $H$ is not an odd minor of $G$, then the chromatic number of $G$ is upper bounded by a number only depending on $H$.

We show that, for any odd minor-closed family, the aforementioned induced bipartite subgraph with large proper conflict-free chromatic number is the only obstruction for having small proper conflict-free chromatic number.

\begin{theorem} \label{thm:pcf_odd_color_intro}
For every positive integer $h$, there exists an integer $c_h$ such that for every graph $G$ with no odd $K_h$-minor, if $\chi_{{\rm pcf}}(Q) \leq k$ for every induced bipartite subgraph $Q$ of $G$, then $\chi_{{\rm pcf}}(G) \leq k+c_h$.
\end{theorem}

In other words, Theorem \ref{thm:pcf_odd_color_intro} implies that every graph with large proper conflict-free chromatic number either contains a large complete graph as an odd minor, or contains an induced bipartite subgraph with large proper conflict-free chormatic number.
It can be viewed as an analog of chi-boundedness result with respect to proper conflict-free chromatic number.

In fact, Theorem \ref{thm:pcf_odd_color_intro} works for a more general setting.
For example, if we want to show $G$ has small odd chromatic number instead of having small proper conflict-free chromatic number, then we can show that asking every induced bipartite subgraph for having small odd chromatic number suffices.
We will state our result in a more general setting, called $S$-achieved coloring, which is a common generalization of conflict-free coloring and odd coloring.

Let $S$ be a set of positive integers.
Let $G$ be a graph.
An {\it $S$-achieved coloring} of $G$ is a coloring such that for every vertex $v$ with $N_G(v) \neq \emptyset$, there exist $s_v \in S$ and a color appearing exactly $s_v$ times in $N_G(v)$.
For a positive integer $t$, we say that $G$ is {\it properly $S$-achieved $t$-colorable} if there exists a proper $S$-achieved $t$-coloring of $G$.
Note that if $1 \in S$, then every graph $H$ is properly $S$-achieved $t$-colorable for some $t \leq |V(H)|$; if $1 \not \in S$, then every graph that has a vertex whose neighborhood is a clique is not properly $S$-achieved $t$-colorable for any integer $t$.
So we are only interested in the case when $1 \in S$.
Moreover, when $S=\{1\}$, the proper $S$-achieved colorings are exactly the proper conflict-free colorings; when $S$ consists of the set of positive odd integers, the proper $S$-achieved colorings are exactly the proper odd colorings.

Similarly, given a list-assignment $L$ of $G$, an {\it $S$-achieved $L$-coloring} of $G$ is an $L$-coloring such that for every vertex $v$ with $N_G(v) \neq \emptyset$, there exist $s_v \in S$ and a color appearing exactly $s_v$ times in $N_G(v)$.
For a positive integer $t$, we say that $G$ is {\it $S$-achieved $t$-choosable} if for every $t$-list-assignment of $G$, there exists a proper $S$-achieved $L$-coloring of $G$.

The following theorem is a generalization of Theorem \ref{thm:pcf_odd_color_intro}.

\begin{theorem} \label{thm:odd_minor_intro}
For every graph $H$, there exists an integer $c_H$ such that the following holds.
Let $G$ be a graph such that $H$ is not an odd minor of $G$.
Let $S$ be a set of positive integers with $1 \in S$.
Let $k$ be a positive integer.
	\begin{enumerate}
		\item If every induced bipartite subgraph of $G$ is properly $S$-achieved $k$-colorable, then $G$ is properly $S$-achieved $(k+c_H)$-colorable.
		\item If every induced bipartite subgraph of $G$ is $S$-achieved $k$-choosable, then $G$ is $S$-achieved $(k+c_H)$-choosable.
	\end{enumerate}
\end{theorem}

This paper is organized as follows.
Results about minor-closed families (Theorems \ref{thm:minor_intro} and \ref{thm:surface_intro}) are proved in Section \ref{sec:minor}.
We develop machinery for conflict-free coloring of clique-sums of graphs in Section \ref{sec:sum} and use it to prove Theorems \ref{thm:topo_intro} and \ref{thm:odd_minor_intro} in Section \ref{sec:topo}.
In Section \ref{sec:layer}, we prove Theorem \ref{thm:ltw_intro}.
Section \ref{sec:layer} does not use anything developed in earlier sections and can be read independently.
Concluding remarks are stated in Section \ref{sec:remarks}.

\section{Minor-closed families} \label{sec:minor}

Let $d$ be a positive real number and let $q$ be a positive integer.
We say that a graph $G$ is {\it $(q,d)$-degenerate} if for every subgraph $H$ of $G$ with $|V(H)| \geq q+1$, $H$ has a vertex of degree at most $d$.
We say a graph class $\F$ is {\it $(q,d)$-degenerate} if every graph in $\F$ is $(q,d)$-degenerate. 
Note that every $d$-degenerate graph class is $(1,d)$-degenerate.

We first prove Theorem \ref{thm:minor_intro}, which is a corollary of the following theorem.
The proof of this theorem uses ideas from \cite{cls}.

\begin{theorem} \label{thm:minor_strong}
Let $d$ be a positive real number and let $q$ be a positive integer.
Let $\F$ be a $(q,d)$-degenerate minor-closed family.
Let $G$ be a graph in $\F$.
Let $L$ be a $\min\{2 \lfloor d \rfloor+1, q\}$-list-assignment of $G$.
Let $H$ be a subgraph of $G$.
Then there exists a proper $L$-coloring $\phi$ of $G$ such that for every vertex $v$ of $G$ with $N_{G-E(H)}(v) \neq \emptyset$, there exists $i$ in the image of $\phi$ such that $|\phi^{-1}(\{i\}) \cap N_{G-E(H)}(v)|=1$.
\end{theorem}

\begin{pf}
We shall prove this theorem by induction on $|V(G)|$.
Let $n=|V(G)|$.
If $n \leq \min\{2 \lfloor d \rfloor+1, q\}$, then there exists an $L$-coloring of $G$ such that no two vertices have the same color, so the theorem holds.
In particular, the theorem holds when $n=1$.
So we may assume that $n \geq \min\{2 \lfloor d \rfloor+1, q\}+1$ and the theorem holds when $n$ is smaller.

Since $\F$ is $(q,d)$-degenerate and $n \geq q+1$, there exists a vertex $v$ of $G$ with degree at most $\lfloor d \rfloor$.
Note that $L|_{V(G)-\{v\}}$ is a $\min\{2 \lfloor d \rfloor+1, q\}$-list-assignment of $G-v$.
Let $H_1=H-v$.
So $H_1$ is a subgraph of $G-v$.
By the induction hypothesis, there exists a proper $L|_{V(G)-\{v\}}$-coloring $\phi_1$ of $G-v$ such that for every vertex $u$ of $G-v$ with $N_{(G-v)-E(H_1)}(u) \neq \emptyset$, there exists $i_1$ in the image of $\phi_1$ such that $|\phi_1^{-1}(\{i_1\}) \cap N_{(G-v)-E(H_1)}(u)|=1$.
If $v$ has no neighbor in $G-E(H)$, then $N_{G-E(H)}(u)=N_{(G-v)-E(H_1)}(u)$ for every $u \in V(G)-\{v\}$, so we can further color $v$ with any color in $L(v)-\{\phi_1(y): y \in N_G(v)\}$ (which is a non-empty set) to extend $\phi_1$ to be a proper $L$-coloring of $G$ satisfying the conclusion of this theorem. 
Hence we may assume that $v$ has at least one neighbor in $G-E(H)$.

Let $u$ be a neighbor of $v$ in $G-E(H)$. 
Let $G_2$ be the graph obtained from $G$ by contracting the edge $uv$ into a new vertex $w$ and deleting resulting parallel edges.
Let $H_2$ be the graph with $V(H_2)=(V(H)-\{u,v\}) \cup \{w\}$ and $E(H_2)=\{e \in E(H): e$ is not incident with $u$ or $v\} \cup \{wz \in E(G_2): z \in V(G)-\{u,v\}, uz \in E(H)\} \cup \{wz \in E(G_2): z \in N_G(v)-\{u\}, uz \not \in E(G)\}$.
Note that $H_2$ is a subgraph of $G_2$, and the edges of $H_2$ are exactly the edges of $H$ remaining in $G_2$ and the new edges obtained by the contraction.
Let $L_2$ be the list-assignment of $G_2$ such that $L_2(w)=L_2(u)$, and $L_2(x)=L(x)$ for every $x \in V(G_2)-\{w\}$.
Note that $L_2$ is a $\min\{2 \lfloor d \rfloor+1, q\}$-list-assignment of $G_2$.

Since $\F$ is a minor-closed family, $G_2 \in \F$.
Since $|V(G_2)|<|V(G)|$, by the induction hypothesis, there exists a proper $L_2$-coloring $\phi_2$ of $G_2$ such that for every vertex $y$ of $G_2$ with $N_{G_2-E(H_2)}(y) \neq \emptyset$, there exists $i_y$ in the image of $\phi_2$ such that $|\phi_2^{-1}(\{i_y\}) \cap N_{G_2-E(H_2)}(y)|=1$.
Let $S=\{\phi_2(y): y \in N_G(v)\} \cup \{i_y: y \in N_G(v), N_{G_2-E(H_2)}(y) \neq \emptyset\}$.
Note that $|S| \leq 2|N_G(v)| \leq 2 \lfloor d \rfloor \leq |L(v)|-1$.
So $L(v)-S \neq \emptyset$.
Then we define an $L$-coloring $\phi$ of $G$ by defining $\phi(x)=\phi_2(x)$ for every $x \in V(G)-\{u,v\}$, defining $\phi(u)=\phi_2(w)$, and defining $\phi(v)$ to be an arbitrary element in $L(v)-S$.
Clearly, $\phi$ is a proper $L$-coloring of $G$.

Note that for every $y \in V(G)-\{u,v\}$, $N_{G_2-E(H_2)}(y)=N_{G-E(H)}(y)-\{v\}$; and $N_{G_2-E(H_2)}(w) \allowbreak =N_{G-E(H)}(u)-\{v\}$.
For convenience, we treat $w$ and $u$ as the same vertex.
Then for every $y \in V(G)-\{v\}$ with $N_{G-E(H)}(y) \neq \emptyset$, either $N_{G_2-E(H_2)}(y) \neq \emptyset$, or $v$ is the unique neighbor of $y$ in $G-E(H)$.
For the former, $|\phi^{-1}(\{i_y\}) \cap N_{G-E(H)}(y)|=|\phi_2^{-1}(\{i_y\}) \cap N_{G_2-E(H_2)}(y)|=1$ since either $y \not \in N_G(v)$ or $\phi(v) \neq i_y \in S$.
For the latter, $|\phi^{-1}(\{\phi(v)\}) \cap N_{G-E(H)}(y)| = |\{v\}|=1$.
Hence, to prove this lemma, it suffices to show that there exists $i$ in the image of $\phi$ such that $|\phi^{-1}(\{i\}) \cap N_{G-E(H)}(v)|=1$.
Recall that $u \in N_{G-E(H)}(v)$.
Since $w$ is adjacent in $G_2$ to all vertices in $N_G(v)-\{u\}$ and $\phi_2$ is a proper $L_2$-coloring of $G_2$, we know $|\phi^{-1}(\{\phi(u)\}) \cap N_{G-E(H)}(v)|=|\{u\}|=1$.
This proves the theorem.
\end{pf}

\bigskip

Now we show corollaries of Theorem \ref{thm:minor_strong}.
Note that the use of the special subgraph $H$ in Theorem \ref{thm:minor_strong} is to make the inductive argument work, and we only need the case $E(H)=\emptyset$ in our applications.
The following corollary is a restatement of Theorem \ref{thm:minor_intro}.

\begin{corollary} \label{cor:minor_cho}
Let $d$ be a nonnegative integer.
Let $\F$ be a $d$-degenerate minor-closed family.
Then every graph $G$ in $\F$ is conflict-free $(2d+1)$-choosable. 
\end{corollary}

\begin{pf}
This corollary trivially holds when $d=0$.
So we may assume $d \geq 1$.
Since $\F$ is $d$-degenerate, $\F$ is $(1,d)$-degenerate.
Let $G$ be a graph in $\F$.
Let $H$ be a subgraph of $G$ with no edge.
So $N_G(v)=N_{G-E(H)}(v)$ for every $v \in V(G)$.
For every $(2d+1)$-list-assignment $L$ of $G$, since $2d+1 \geq \min\{2d+1,1\}$, by Theorem \ref{thm:minor_strong}, there exists a proper conflict-free $L$-coloring of $G$.
\end{pf}

\bigskip

The following corollary is a restatement of Theorem \ref{thm:surface_intro}.

\begin{corollary}
Let $\Sigma$ be a surface with Euler genus $\rho$.
If $\rho \in \{0,1\}$, then every graph embeddable in $\Sigma$ is conflict-free $11$-choosable.
If $\rho \geq 2$, then every graph embeddable in $\Sigma$ is conflict-free $\frac{13+\sqrt{73+48\rho}}{2}$-choosable.
\end{corollary}

\begin{pf}
Let $g$ be the function such that $g(x)=6-\frac{12-6\rho}{x}$ for every positive integer $x$.
Let $\F$ be the class of graphs embeddable in $\Sigma$.
By Euler's formula, every $n$-vertex graph in $\F$ is $g(n)$-degenerate. 

If $\rho \in \{0,1\}$, then $\F$ is 5-degenerate, so every graph $G$ in $\F$ is $11$-conflict-free choosable by Corollary \ref{cor:minor_cho}.
Hence we may assume $\rho \geq 2$.

Let $k=\frac{13+\sqrt{73+48\rho}}{2}$.
If $G$ is a subgraph of a graph in $\F$ on $n$ vertices, for some integer $n$ with $n \geq \lfloor k \rfloor+1$, then $n \geq k$, so $g(n) = 6+\frac{6\rho-12}{n} \leq 6+\frac{6\rho-12}{k} \leq \frac{k-1}{2}$.
Hence $\F$ is $(\lfloor k \rfloor,\frac{k-1}{2})$-degenerate.
Since $k \geq \min\{2 \lfloor \frac{k-1}{2} \rfloor+1, \lfloor k \rfloor\}$, this corollary follows from Theorem \ref{thm:minor_strong} with $H$ the null graph.
\end{pf}

\section{$S$-achieved coloring and clique-sum} \label{sec:sum} 

In this section we study conflict-free coloring for graphs with a given tree-decomposition.
It is preparation for proving Theorems \ref{thm:topo_intro} and \ref{thm:odd_minor_intro}.

We first define clique-sums of graphs.
Let $G_1,G_2$ be graphs.
For each $i \in [2]$, let $Q_i$ be a clique in $G_i$ with $|Q_1|=|Q_2|$.
Let $\iota$ be a bijection from $Q_1$ to $Q_2$.
Then a {\it $(Q_1,Q_2,\iota)$-sum of $G_1$ and $G_2$} is a graph obtained from a disjoint union of $G_1$ and $G_2$ by first for each $x \in Q_1$, identifying $x$ and $\iota(x)$ into a new vertex $v_x$, and then deleting any number of edges whose both ends are in $\{v_x: x \in Q_1\}$.
Furthermore, if $L_1$ is a list-assignment of $G_1$, and $L_2$ is a list-assignment of $G_2$ such that $L_1(x)=L_2(\iota(x))$ for every $x \in Q_1$, then {\it the list-assignment of $G$ obtained by a $(Q_1,Q_2,\iota)$-sum} is the list-assignment $L$ of $G$ such that $L(v)=L_1(v)$ if $v \in V(G_1)-Q_1$, $L(v)=L_2(v)$ if $v \in V(G_2)-Q_2$, and $L(v_x)=L_1(x)=L_2(\iota(x))$ for every $x \in Q_1$.
For a nonnegative integer $t$, an {\it $(\leq t)$-sum} of $G_1$ and $G_2$ is a $(W_1,W_2,\xi)$-sum of $G_1$ and $G_2$ for some clique $W_1$ in $G_1$, some clique $W_2$ in $G_2$ with $|W_1|=|W_2| \leq t$, and a bijection $\xi$ from $W_1$ to $W_2$.

In order to handle clique-sums of graphs, we consider the following property stronger than $S$-achieved colorability.

Let $S$ be a set of positive integers.
Let $G$ be a graph.
Let $\C$ be a collection of subsets of $V(G)$.
We say that a subgraph $H$ of $G$ is {\it $\C$-compatible} if for every edge of $H$, there exists a member of $\C$ containing both ends of this edge.
For a nonnegative integer $t$ and a list-assignment $L$ of $G$, we say that $G$ is {\it $(S,\C,L,t)$-extendable} if for every clique $W$ of $G$ with size at most $t$, for every $\C$-compatible subgraph $H$ of $G$, for every proper $L|_W$-coloring $\phi_W$ of $G[W]$, for every function $f$ with domain $W$, and for every function $g: W \rightarrow \{0,1\}$, there exists a proper $L$-coloring $\phi$ of $G$ such that
	\begin{itemize}
		\item $\phi(v)=\phi_W(v)$ for every $v \in W$, and
		\item for every $v \in V(G)$,
			\begin{itemize}
				\item if $v \in W$ and $g(v)=0$, then there exists no vertex $u \in N_G(v)-W$ such that $\phi(u)=f(v)$, and
				\item if either $v \not \in W$, or $v \in W$ and $g(v)=1$, then either $N_{G-E(H)}(v)=\emptyset$, or there exists $i_v$ such that $|\phi^{-1}(\{i_v\}) \cap N_{G-E(H)}(v)| \in S$.
			\end{itemize}
	\end{itemize}
Note that if $G$ is $(S,\C,L,t)$-extendable, then $G$ has a proper $S$-achieved $L$-coloring, by taking $W=\emptyset$ and $H$ the null graph.

\begin{lemma} \label{lemma:clique_sum}
Let $S$ be a set of positive integers.
Let $t$ be a nonnegative integer.
Let $G_1$ and $G_2$ be graphs.
For each $i \in [2]$, let $Q_i$ be a clique in $G_i$ with size $t$, and let $\C_i$ be a collection of subsets of $V(G_i)$ with $Q_i \in \C_i$.
Let $\iota$ be a bijection from $Q_1$ to $Q_2$.
For each $i \in \{1,2\}$, let $L_i$ be a list-assignment of $G_i$ such that $L_1(x)=L_2(\iota(x))$ for every $x \in Q_1$.
Let $G$ be a $(Q_1,Q_2,\iota)$-sum of $G_1$ and $G_2$, and let $L$ be the list-assignment obtained by a $(Q_1,Q_2,\iota)$-sum. 
If $G_1$ is $(S,\C_1,L_1,t_1)$-extendable and $G_2$ is $(S,\C_2,L_2,t_2)$-extendable for some integers $t_1 \geq t$ and $t_2 \geq t$, then $G$ is $(S,\C_1 \cup \C_2,L,\min\{t_1,t_2\})$-extendable.
\end{lemma}

\begin{pf}
Let $G'$ be the graph obtained from a disjoint union of $G_1$ and $G_2$ by, for each vertex $x \in Q_1$, identifying $x$ and $\iota(x)$ into a vertex $v_x$.
Let $Q=\{v_x: x \in Q_1\}$.
So $Q$ is a clique in $G'$, and $G$ is a graph obtained from $G'$ by deleting some edges whose both ends are in $Q$.
For convenience, we do not distinguish $v_x$, $x$ and $\iota(x)$.
That is, we treat $v_x=x=\iota(x)$.
This implies that $Q=Q_1=Q_2$ is a clique in all of $G_1$, $G_2$ and $G'$.

Let $W$ be a clique of $G$ with size at most $\min\{t_1,t_2\}$.
Let $H$ be a $(\C_1 \cup \C_2)$-compatible subgraph of $G$.
Let $\phi_W$ be a proper $L|_W$-coloring of $G[W]$.
Let $f$ be a function with domain $W$.
Let $g: W \rightarrow \{0,1\}$ be a function.
To show that $G$ is $(S,L,\min\{t_1,t_2\})$-extendable, it suffices to show that there exists a proper $L$-coloring $\phi$ of $G$ such that the following hold.
	\begin{itemize}
		\item[(i)] $\phi(v)=\phi_W(v)$ for every $v \in W$.
		\item[(ii)] For every $v \in V(G)$,
			\begin{itemize}
				\item[(ii-1)] if $v \in W$ and $g(v)=0$, then there exists no vertex $u \in N_G(v)-W$ such that $\phi(u)=f(v)$, and
				\item[(ii-2)] if either $v \not \in W$, or $v \in W$ and $g(v)=1$, then either $N_{G-E(H)}(v)=\emptyset$, or there exists $i_v$ such that $|\phi^{-1}(\{i_v\}) \cap N_{G-E(H)}(v)| \in S$.
			\end{itemize}
	\end{itemize}

Since $W$ is a clique of $G \subseteq G'$, either $W \subseteq V(G_1)$ or $W \subseteq V(G_2)$.
By symmetry, we may assume $W \subseteq V(G_1)$.
So $W$ is a clique in both $G$ and $G_1$.
Hence $\phi_W$ is a proper $L_1|_W$-coloring of $G_1[W]$.
Let $H_1$ be the subgraph of $G_1$ with $V(H_1)=V(H) \cap V(G_1)$ and $E(H_1)=\{uv \in E(H): u,v \in V(G_1)\} \cup \{xy: x,y \in Q_1, v_xv_y \not \in E(G)\}$.
Since $Q_1 \in \C_1$, $H_1$ is a $\C_1$-compatible subgraph of $G_1$, and $G_1-E(H_1)=(G-E(H))[V(G_1)]$.
Since $G_1$ is $(S,\C_1,L_1,t_1)$-extendable with $t_1 \geq |W|$, there exists a proper $L_1$-coloring $\phi_1$ of $G_1$ such that $\phi_1|_W=\phi_W$ and for every $v \in V(G_1)$, 
	\begin{itemize}
		\item if $v \in W$ and $g(v)=0$, then there exists no vertex $u \in N_{G_1}(v)-W$ such that $\phi_1(u)=f(v)$, and
		\item if either $v \not \in W$, or $v \in W$ and $g(v)=1$, then either $N_{G_1-E(H_1)}(v)=\emptyset$, or there exists $i_v$ such that $|\phi_1^{-1}(\{i_v\}) \cap N_{G_1-E(H_1)}(v)| \in S$.
	\end{itemize}

Let $H_2$ be the subgraph of $G_2$ with $V(H_2)=V(H) \cap V(G_2)$ and $E(H_2)=\{uv \in E(H): u,v \in V(G_2)\} \cup \{\iota(x)\iota(y): x,y \in Q_2,v_xv_y \not \in E(G)\}$.
So $G_2-E(H_2)=(G-E(H))[V(G_2)]$.
Since $Q_2 \in \C_2$, $H_2$ is $\C_2$-compatible.
Note that $Q=Q_1=Q_2$ is a clique in both $G_1$ and $G_2$ with size $t \leq t_2$, and $\phi_1|_Q$ is a proper $L_2$-coloring of $G_2[Q_2]$.
Define $g_2: Q_2 \rightarrow \{0,1\}$ to be the function and $f_2$ to be the function with domain $Q_2$ such that for every $x \in Q_2$,
	\begin{itemize}
		\item if $x \in Q_2 \cap W$, and either $N_{G_1-E(H_1)}(x)=\emptyset$ or $g(x)=0$, then $g_2(x)=g(x)$ and $f_2(x)=f(x)$,
		\item if $x \in Q_2 \cap W$, $N_{G_1-E(H_1)}(x) \neq \emptyset$ and $g(x)=1$, then $g_2(x)=0$ and $f_2(x)=i_x$,
		\item if $x \in Q_2 \setminus W$ and $N_{G_1-E(H_1)}(x)=\emptyset$, then $g_2(x)=1$ and $f_2(x)=0$, and
		\item if $x \in Q_2 \setminus W$ and $N_{G_1-E(H_1)}(x) \neq \emptyset$, then $g_2(x)=0$ and $f_2(x)=i_x$.
	\end{itemize}
Since $G_2$ is $(S,\C_2,L_2,t_2)$-extendable, there exists a proper $L_2$-coloring $\phi_2$ of $G_2$ such that $\phi_2|_{Q_2}=\phi_1|_{Q_1}$ and for every $v \in V(G_2)$, 
	\begin{itemize}
		\item if $v \in Q_2$ and $g_2(v)=0$, then there exists no vertex $u \in N_{G_2}(v)-Q_2$ such that $\phi_2(u)=f_2(v)$, and
		\item if either $v \not \in Q_2$, or $v \in Q_2$ and $g_2(v)=1$, then either $N_{G_2-E(H_2)}(v)=\emptyset$, or there exists $j_v$ such that $|\phi_2^{-1}(\{j_v\}) \cap N_{G_2-E(H_2)}(v)| \in S$.
	\end{itemize}

Since $\phi_2|_{Q_2}=\phi_1|_{Q_1}$, the function $\phi$ with domain $V(G)$ with $\phi|_{V(G_1)}=\phi_1$ and $\phi|_{V(G_2)}=\phi_2$ is a proper $L$-coloring of $G$ such that $\phi|_W=\phi_1|_W=\phi_W$.
So (i) holds.

Now we show that (ii) holds.
We first prove (ii-1).
Let $a \in W$ with $g(a)=0$.
By the property of $\phi_1$, there exists no vertex $u \in N_{G_1}(a)-W$ such that $\phi_1(u)=f(a)$.
If $a \not \in Q$, then $N_G(a)-W=N_{G_1}(a)-W$, so there exists no vertex $u \in N_{G}(a)-W$ such that $\phi(u)=f(a)$.
So we may assume $a \in Q \cap W=Q_2 \cap W$.
Since $g(a)=0$, we have $g_2(a)=0$ and $f_2(a)=f(a)$.
By the property of $\phi_2$, there exists no vertex $u \in N_{G_2}(a)-Q_2$ such that $\phi_2(u)=f_2(a)=f(a)$.
Since $N_G(a) \subseteq N_{G_1}(a) \cup N_{G_2}(a)$, there exists no vertex $u \in N_{G}(a)-W$ such that $\phi(u)=f(a)$.
This proves (ii-1).

Now we prove (ii-2).
Let $b$ be a vertex of $G$ such that either $b \not \in W$, or $b \in W$ and $g(b)=1$.
To prove (ii), it suffices to show that either $N_{G-E(H)}(b)=\emptyset$, or there exists $k_b$ such that $|\phi^{-1}(\{k_b\}) \cap N_{G-E(H)}(b)| \in S$.
We may assume $N_{G-E(H)}(b) \neq \emptyset$, for otherwise we are done.

We first assume $b \in V(G_2)-V(G_1)$.
In particular, $b \not \in Q$.
By the property of $\phi_2$, either $N_{G_2-E(H_2)}(b)=\emptyset$, or $j_b$ exists and $|\phi_2^{-1}(\{j_b\}) \cap N_{G_2-E(H_2)}(b)| \in S$.
Since $b \in V(G_2)-V(G_1)$, $N_{G_2-E(H_2)}(b)=N_{G-E(H)}(b) \neq \emptyset$.
So we are done by choosing $k_b=j_b$.

Hence we may assume $b \in V(G_1)$.
Since either $b \not \in W$, or $b \in W$ and $g(b)=1$, by the property of $\phi_1$, we know that either $N_{G_1-E(H_1)}(b)=\emptyset$, or $i_b$ exists and $|\phi_1^{-1}(\{i_b\}) \cap N_{G_1-E(H_1)}(b)| \in S$.
Since $(G-E(H))[V(G_1)]=G_1-E(H_1)$, we are done if $b \not \in Q_1$.
So we may assume $b \in Q_1 =Q_2$.
Hence $N_{G-E(H)}(b)=N_{G_1-E(H_1)}(b) \cup N_{G_2-E(H_2)}(b)$. 

If $N_{G_1-E(H_1)}(b)=\emptyset$, then $g_2(b)=1$, so the property of $\phi_2$ implies that either $N_{G_2-E(H_2)}(b)=\emptyset$, or $j_b$ exists and $|\phi_2^{-1}(\{j_b\}) \cap N_{G_2-E(H_2)}(b)| \in S$.
In this case, $N_{G_2-E(H_2)}(b)=N_{G_1-E(H_1)}(b) \cup N_{G_2-E(H_2)}(b) = N_{G-E(H)}(b) \neq \emptyset$, so the latter holds, and $|\phi^{-1}(\{j_b\}) \cap N_{G-E(H)}(b)| = |\phi_2^{-1}(\{j_b\}) \cap N_{G_2-E(H_2)}(b)| \in S$.

So we may assume $N_{G_1-E(H_1)}(b) \neq \emptyset$.
Recall that either $b \not \in W$, or $b \in W$ and $g(b)=1$.
Hence the property of $\phi_1$ implies that $i_b$ exists and $|\phi_1^{-1}(\{i_b\}) \cap N_{G_1-E(H_1)}(b)| \in S$.
If $b \not \in W$, then $b \in Q_2 \setminus W$ and $N_{G_1-E(H_1)}(b) \neq \emptyset$, so $g_2(b)=0$ and $f_2(b)=i_b$; if $b \in W$ and $g(b)=1$, then $b \in Q_2 \cap W$, $N_{G_1-E(H_1)}(b) \neq \emptyset$ and $g(b)=1$, so  $g_2(b)=0$ and $f_2(b)=i_b$.
In either case, $g_2(b)=0$ and $f_2(b)=i_b$.
So the property of $\phi_2$ implies that $|\phi_2^{-1}(\{i_b\}) \cap N_{G_2-E(H_2)}(b)-Q_2| \leq |\phi_2^{-1}(\{i_b\}) \cap N_{G_2}(b)-Q_2|=0$.
Therefore, $|\phi^{-1}(\{i_b\}) \cap N_{G-E(H)}(b)| = |\phi_1^{-1}(\{i_b\}) \cap N_{G_1-E(H_1)}(b)|+|\phi_2^{-1}(\{i_b\}) \cap N_{G_2-E(H_2)}(b)-Q_2| = |\phi_1^{-1}(\{i_b\}) \cap N_{G_1-E(H_1)}(b)| \in S$.
This proves the lemma.
\end{pf}

\bigskip

A \emph{tree-decomposition} of a graph $G$ is a sequence $(B_x:x\in V(T))$ of subsets of $V(G)$ called \emph{bags} that are indexed by the nodes of a tree $T$ such that
\begin{itemize}
	\item for each vertex $v$ of $G$, $T[\{x\in V(T):v\in B_x\}]$ is a connected graph with at least one vertex; and
	\item for each edge $vw$ of $G$, there exists a node $x$ of $T$ such that $\{v,w\}\subseteq B_x$.
\end{itemize}
The \emph{width} of a tree-decomposition is equal to the size of its largest bag, minus $1$.  
The \emph{treewidth} of $G$ is the minimum width of any tree-decomposition of $G$.

Let $G$ be a graph.
Let $\T=(B_x: x \in V(T))$ be a tree-decomposition of $G$.
The {\it adhesion} of $\T$ is $\max_{xy \in E(T)}|B_x \cap B_y|$.
Let $x \in V(T)$.
Then the {\it torso} at $x$ (with respect to $\T$) is the graph obtained from $G[B_x]$ by for each edge $xx' \in E(T)$, adding edges such that $B_x \cap B_{x'}$ is a clique.
The {\it frame} at $x$ (with respect to $\T$) is the collection $\{B_x \cap B_{x'}: xx' \in E(T)\}$.
It is simple to see that $G$ is an $(\leq \xi)$-sum of the torsos of $\T$, where $\xi$ is the adhesion of $\T$, and the cliques involved in the clique-sums are the members of the frames at the nodes for the corresponding torsos.
For a set of positive integers $S$ and a list-assignment $L$ of $G$, we say that $\T$ is {\it $L$-extendable} if for every $x \in V(T)$, $R_x$ is $(S,\C_x,L|_{V(R_x)},\xi)$-extendable, where $R_x$ is the torso at $x$, $\C_x$ is the frame at $x$, and $\xi$ is the adhesion of $\T$.

\begin{lemma} \label{lemma:torsos_ext}
Let $S$ be a set of positive integers.
Let $\xi$ be a nonnegative integer.
Let $G$ be a graph.
Let $\T$ be a tree-decomposition of $G$ with adhesion at most $\xi$. 
Let $L$ be a list-assignment of $G$.
If $\T$ is $L$-extendable, then $G$ is $(S,\emptyset,L,\xi)$-extendable.
\end{lemma}

\begin{pf}
Let $t_1,t_2,...,t_{|V(T)|}$ be a breadth-first-search ordering of $V(T)$.
For every $1 \leq i \leq |V(T)|$, let $R_i$ be the torso at $t_i$ and $\C_i$ the frame at $t_i$.
Let $G^+$ be the graph obtained from $G$ by adding edges such that every member of $\bigcup_{i=1}^{|V(T)|}\C_i$ is a clique.
Since $\T$ is $L$-extendable, each $R_i$ is $(S,\C_i,L|_{V(R_i)},\xi)$-extendable.
By induction on $i$, Lemma \ref{lemma:clique_sum} implies that for every $i$, $G^+[\bigcup_{j=1}^iV(R_i)]$ is $(S,\bigcup_{j=1}^i\C_i,L|_{\bigcup_{j=1}^iV(R_i)},\xi)$-extendable.
So $G^+$ is $(S,\bigcup_{j=1}^{|V(T)|}\C_i,L,\xi)$-extendable.
Let $H$ be the graph with $V(H)=V(G)$ and $E(H)=E(G^+)-E(G)$.
Then $H$ is a $\bigcup_{j=1}^{|V(T)|}\C_i$-compatible subgraph of $G^+$, and $G^+-E(H)=G$.
Hence $G$ is $(S,\emptyset,L,\xi)$-extendable.
\end{pf}

\bigskip

Let $G$ be a graph and $L$ a list-assignment of $G$.
For a set $Z$, we say that a list-assignment $L'$ of $G$ is {\it $(L,Z)$-removed} if for every $v \in V(G)$, $L'(v) \subseteq L(v)$ and $L(v)-L'(v) \subseteq Z$.

\begin{lemma} \label{lemma:cho_to_ext}
Let $S$ be a set of positive integers with $1 \in S$.
Let $G$ be a graph.
Let $\C$ be a collection of subsets of $V(G)$.
Let $\xi$ be a nonnegative integer.
Let $L$ be a list-assignment of $G$. 
If $G-E(H)$ has a proper $S$-achieved $L_H$-coloring for every $\C$-compatible subgraph $H$ of $G$ and every list-assignment $L_H$ of $G$ that is $(L,Z)$-removed for some set $Z$ with $|Z| \leq 2\xi$, then $G$ is $(S,\C,L,\xi)$-extendable. 
\end{lemma}

\begin{pf}
Let $W$ be a clique of size at most $\xi$.
Let $H$ be a $\C$-compatible subgraph of $G$.
Let $\phi_W$ be a proper $L|_W$-coloring of $G[W]$.
Let $f$ be a function with domain $W$.
Let $g: W \rightarrow \{0,1\}$ be a function.
To prove this lemma, if suffices to show that there exists a proper $L$-coloring of $G$ satisfying the conditions in the definition of being $(S,\C,L,\xi)$-extendable.

Let $G'=G-E(H)$.
Let $Z=\{\phi_W(w),f(w): w \in W\}$.
So $|Z| \leq 2\xi$.
For every $v \in V(G')$, let $L'(v)=L(v)-Z$.
So $L'$ is an $(L,Z)$-removed list-assignment of $G'$.
Since $H$ is $\C$-compatible, $G'=G-E(H)$ has a proper $S$-achieved $L'$-coloring $\phi'$.
Let $\phi$ be the $L$-coloring of $G$ such that $\phi|_W=\phi_W$ and $\phi|_{V(G)-W}=\phi'|_{V(G')-W}$.
By the definition of $L'$, the image of $\phi_W$ and the image of $\phi'$ are disjoint, so $\phi$ is a proper $L$-coloring with $\phi|_W=\phi_W$.

Let $v \in V(G)$.
If $v \in W$, then $f(v) \in Z$, so there exists no vertex $u \in N_G(v)-W$ such that $\phi(u)=f(v)$.
Since $\phi'$ is a proper $S$-achieved $L'$-coloring of $G'=G-E(H)$, either $N_{G-E(H)}(v)=N_{G'}(v)=\emptyset$, or there exists $i_v$ such that $|{\phi'}^{-1}(\{i_v\}) \cap N_{G-E(H)}(v)| \in S$.
So if either $v \in W$ and $g(v)=0$, or $N_{G-E(H)}(v)=\emptyset$, then $v$ satisfies the condition mentioned in the definition of being $(S,\C,L,\xi)$-extendable. 
Hence we may assume that $N_{G-E(H)}(v) \neq \emptyset$ and either $v \not \in W$, or $v \in W$ and $g(v)=1$.
So $i_v$ exists.
Since $Z$ is disjoint from the image of $\phi'$, $i_v$ is not in the image of $\phi_W$.
Hence if $N_{G-E(H)}(v) \cap W=\emptyset$, then $|{\phi}^{-1}(\{i_v\}) \cap N_{G-E(H)}(v)| = |{\phi'}^{-1}(\{i_v\}) \cap N_{G-E(H)}(v)| \in S$, so we are done.
So we may assume that there exists $u \in N_{G-E(H)}(v) \cap W$. 
Since $\phi_W$ is a proper $L|_W$-coloring of the subgraph induced by the clique $W$, and $\phi(u)=\phi_W(u)$ is not in the image of $\phi'$, we know $|\phi^{-1}(\{\phi(u)\}) \cap N_{G-E(H)}(v)|=|\{u\}|=1 \in S$.
This proves the lemma. 
\end{pf}

\bigskip

A {\it design} is a pair $(G,\C)$, where $G$ is a graph and $\C$ is a collection of subsets of $V(G)$.
The {\it rank} of a design $(G,\C)$ is $\max_{C \in \C}|C|$.
We say that two designs $(G_1,\C_1)$ and $(G_2,\C_2)$ are {\it isomorphic} if there exist an isomorphism $\iota$ from $G_1$ to $G_2$ and a bijection $\eta$ from $\C_1$ to $\C_2$ such that for every $Y \in \C_1$, $\eta(Y)=\{\iota(y): y \in Y\}$.

Let $\D$ be a set of designs.
Let $\T$ be a tree-decomposition of a graph $G$.
We say that {\it $\T$ is over $\D$} if for every $x \in V(T)$, $(R_x,\C_x)$ is isomorphic to a member of $\D$, where $R_x$ is the torso at $x$ and $\C_x$ is the frame at $x$.

\begin{lemma} \label{lemma:design_torsos}
Let $S$ be a set of positive integers with $1 \in S$.
Let $\D$ be a set of designs.
Let $G$ be a graph.
Let $t$ and $\xi$ be nonnegative integers.
Let $\T$ be a tree-decomposition of $G$ over $\D$ with adhesion at most $\xi$.
If for every member $(Q,\C)$ of $\D$, $Q-E(H)$ is $S$-achieved $t$-choosable (and properly $S$-achieved $t$-colorable, respectively) for every $\C$-compatible subgraph $H$ of $Q$, then $G$ is $S$-achieved $(t+2\xi)$-choosable (and properly $S$-achieved $(t+2\xi)$-colorable, respectively).
\end{lemma}

\begin{pf}
For each $(Q,\C) \in \D$, since $Q-E(H)$ is $S$-achieved $t$-choosable (and properly $S$-achieved $t$-colorable, respectively) for every $\C$-compatible subgraph $H$ of $Q$, by Lemma \ref{lemma:cho_to_ext}, $Q$ is $(S,\C,L,\xi)$-extendable for every $(t+2\xi)$-list-assignment $L$ of $Q$ (and for the $(t+2\xi)$-list-assignment $L$ of $Q$ with $L(v)=\{1,2,...,t+2\xi\}$ for $v \in V(Q)$, respectively).

Let $L$ be a $(t+2\xi)$-list-assignment of $G$ (and the list-assignment $L$ with $L(v)=\{1,2,...,t+2\xi\}$ for every $v \in V(G)$, respectively).
Let $\D'=\{(R_x,\C_x): x \in V(T)\}$, where $R_x$ is the torso at $x$ and $\C_x$ is the frame at $x$.
Then $\T$ is a tree-decomposition of $G$ over $\D'$ with adhesion at most $\xi$ such that for every $(Q,\C) \in \D'$, $Q$ is $(S,\C,L|_{V(Q)},\xi)$-extendable.
Hence $\T$ is $L$-extendable.
So $G$ is $(S,\emptyset,L,\xi)$-extendable by Lemma \ref{lemma:torsos_ext}.
Hence $G$ has a proper $S$-achieved $L$-coloring.
This proves the lemma.
\end{pf}

\begin{lemma} \label{lemma:torsos}
Let $S$ be a set of positive integers with $1 \in S$.
Let $G$ be a graph.
Let $t$ and $\xi$ be nonnegative integers.
Let $\T$ be a tree-decomposition of $G$ with adhesion at most $\xi$.
If for every torso $R$ with respect to $\T$, every subgraph of $R$ is $S$-achieved $t$-choosable, then $G$ is $S$-achieved $(t+2\xi)$-choosable.
\end{lemma}

\begin{pf}
Let $\D = \{(R_x,\C_x): x \in V(T)\}$, where $R_x$ is the torso at $x$ and $\C_x$ is the frame at $x$.
Then this lemma follows from Lemma \ref{lemma:design_torsos}.
\end{pf}

\section {Topological minors and odd minors} \label{sec:topo}

We will prove Theorems \ref{thm:topo_intro} and \ref{thm:odd_minor_intro} by using the machinery developed in the previous section.
We first prove some simple lemmas.

\begin{lemma} \label{lemma:deletion}
Let $S$ be a set of integers with $1 \in S$.
Let $G$ be a graph.
Let $k$ and $c$ be nonnegative integers.
Let $Y$ be a subset of $V(G)$ with $|Y| \leq k$.
If $G-X$ is $S$-achieved $c$-choosable (and properly $S$-achieved $c$-colorable, respectively) for every set $X$ with $X \supseteq Y$ and $|X| \leq 2|Y|$, then $G$ is $S$-achieved $(c+2k)$-choosable (and properly $S$-achieved $(c+2k)$-colorable, respectively).
\end{lemma}

\begin{pf}
Let $L$ be a $(c+2k)$-list-assignment of $G$ (and the $(c+2k)$-list-assignment of $G$ with $L(v)=\{1,2,...,c+2k\}$ for every $v \in V(G)$, respectively).
We shall define a proper $S$-achieved $L$-coloring of $G$.
Without loss of generality, we may assume that $G$ has no isolated vertices.

For every $y \in Y$, since $G$ has no isolated vertices, there exists a neighbor $f(y)$ of $y$ in $G$.
Let $X=\{y,f(y): y \in Y\}$.
So $|X| \leq 2|Y| \leq 2k$.
Hence there exists a proper $L|_X$-coloring $\phi_X$ of $G[X]$ such that all vertices in $X$ use different colors.
Since $1 \in S$, $\phi_X$ is a proper $S$-achieved coloring of $G[X]$.

Let $Z=\{\phi_X(x): x \in X\}$.
So $|Z| = |X| \leq 2k$.
For every $v \in V(G)-X$, let $L'(v)=L(v)-Z$.
Let $G'=G-X$.
Note that for every $v \in V(G')$, $|L'(v)| \geq |L(v)|-|Z| \geq c$.
Since $|X| \leq 2|Y|$ and $X \supseteq Y$, $G'$ has a proper $S$-achieved $L'$-coloring $\phi'$ by assumption.
Let $\phi$ be the $L$-coloring of $G$ such that $\phi|_X=\phi_X$ and $\phi|_{V(G)-X}=\phi'$.
Since the image of $\phi_X$ is $Z$, which is disjoint from the image of $\phi'$, $\phi$ is a proper $L$-coloring of $G$.

Since $G[X]$ has no isolated vertices, for every $v \in X$, there exists $i_v \in Z$ such that $|\phi^{-1}(\{i_v\}) \cap N_G(v)| = |\phi^{-1}(\{i_v\}) \cap N_{G[X]}(v)| = 1 \in S$.
Let $v \in V(G)-X$.
If $v$ is adjacent in $G$ to some vertex $x$ in $X$, then $|\phi^{-1}(\{\phi(x)\}) \cap N_G(v)| = |\{x\}| = 1 \in S$. 
If $v$ is not adjacent in $G$ to any vertex in $X$, then since $G$ has no isolated vertex, $v$ is not an isolated vertex in $G'$, so there exists $i_v \not \in Z$ such that $|\phi^{-1}(\{i_v\}) \cap N_G(v)| = |{\phi'}^{-1}(\{i_v\}) \cap N_{G'}(v)| \in S$. 
Therefore, $\phi$ is a proper $S$-achieved $L$-coloring. 
This proves the lemma.
\end{pf}

\begin{lemma} \label{lemma:near_bdd_max_deg}
Let $d \geq 1$ be a real number.
If $G$ is a graph that has at most $d$ vertices with degree at least $d$, then $G$ is conflict-free $(4d-1)$-choosable. 
\end{lemma}

\begin{pf}
Let $Y$ be the set of vertices in $G$ with degree at least $d$.
So $|Y| \leq d$.
For every set $X$ with $X \supseteq Y$, $G-X$ has maximum degree at most $d-1$, so by \cite[Proposition 3]{cl}, $G-X$ is conflict-free $(2d-1)$-choosable.\footnote{Proposition 3 in \cite{cl} is stated for coloring instead of list-coloring. But it is easy to see that the simple proof of \cite[Proposition 3]{cl} works for list-coloring. Or one can also avoid any known result in the literature here by just proving a weaker form of this lemma by increasing the number $4d-1$ stated in the statement of this lemma to $d^2+2$. This weaker form is still strong enough to prove Theorem \ref{thm:topo_intro}, and it immediately follows from the trivial fact that $(G-X)^2$ is $((d-1)^2+1)$-choosable.}
Then this lemma follows from Lemma \ref{lemma:deletion}.
\end{pf}

\bigskip

Now we are ready to prove Theorem \ref{thm:topo_intro}.
The following is a restatement.

\begin{theorem} \label{thm:topo}
For every graph $H$, there exists a real number $c_H$ such that every graph that does not contain a subdivision of $H$ is conflict-free $c_H$-choosable.
\end{theorem}

\begin{pf}
Let $H$ be a graph.
By \cite[Theorem 4.1]{gm}, there exists a positive integer $c_1$ (only depending on $H$) such that every graph that does not contain a subdivision of $H$ has a tree-decomposition $\T=(B_x: x \in V(T))$ with adhesion at most $c_1$ such that for every $x \in V(T)$, the torso at $x$ either contains at most $c_1$ vertices with degree at least $c_1$ or does not contain $K_{c_1}$ as a minor.
By \cite{m_minor}, there exists a real number $c_2$ (only depending on $c_1$ and hence only depending on $H$) such that every graph with no $K_{c_1}$-minor is $c_2$-degenerate. 
Define $c_H=6c_1+2c_2-1$.

Let $\F_1$ be the class of graphs that have at most $c_1$ vertices with degree at least $c_1$.
Let $\F_2$ be the class of $K_{c_1}$-minor free graphs. 
Let $\F=\F_1 \cup \F_2$.
Note that every subgraph of a graph in $\F$ is in $\F$.
By Lemma \ref{lemma:near_bdd_max_deg}, every graph in $\F_1$ is conflict-free $(4c_1-1)$-choosable.
By Corollary \ref{cor:minor_cho}, every graph in $\F_2$ is conflict-free $(2c_2+1)$-choosable.
Hence every graph in $\F$ is $\{1\}$-achieved $\max\{4c_1-1,2c_2+1\}$-choosable and hence $\{1\}$-achieved $(4c_1+2c_2-1)$-choosable.

Let $G$ be a graph that does not contain a subdivision of $H$.
Then $G$ has a tree-decomposition with adhesion at most $c_1$ such that for every torso $R$, every subgraph of $R$ is in $\F$ and hence is $\{1\}$-achieved $(4c_1+2c_2-1)$-choosable. 
By Lemma \ref{lemma:torsos}, $G$ is $\{1\}$-achieved $(6c_1+2c_2-1)$-choosable. 
Therefore, $G$ is conflict-free $c_H$-choosable.
\end{pf}

\bigskip

Now we prove results for odd minors.
We need the following structure theorem for odd minors (Theorem \ref{thm:odd_structure}).
This theorem is a simple combination of known results in the literature.
It is likely folklore but seems not formally written in the literature.
So we provide a proof here for completeness.

A {\it torso} of a design $(G,\C)$ is the graph obtained from $G$ by adding edges such that for every added edge, some member of $\C$ contains the both ends of this edge.

\begin{theorem} \label{thm:odd_structure}
For every graph $H$, there exist positive integers $r,\xi$ such that if $H$ is not an odd minor of a graph $G$, then $G$ has a tree-decomposition $\T=(B_x: x \in V(T))$ for some tree $T$ such that for every $x \in V(T)$, either $K_r$ is not a minor of the torso at $x$, or there exists $Z \subseteq V(R_x)$ with $|Z| \leq \xi$ such that $R_x-Z$ is an induced bipartite subgraph of $G$ and $|M-Z| \leq 1$ for every member $M$ of the frame at $x$, where $R_x$ is the torso at $x$.
\end{theorem}

\begin{pf}
This proof involves many notions that will be not used in the rest of the paper, so we do not include their formal definition here.
Among those notions, ``$K_h$-expansion'' is defined in \cite{ggrsv} and all other undefined notions are the ones defined in \cite{rs}.

Let $h=|V(H)|$.
By \cite[Theorem 13]{ggrsv} (taking $\ell=h$), for every graph $G'$ with no odd $K_h$-minor, there exists a positive integer $t$ (only depending on $h$) with $t \geq 12h$ such that if $G'$ contains a $K_t$-expansion $\eta$, then there exist $X \subseteq V(G')$ with $|X| <8h$ and an induced bipartite subgraph $U$ of $G'-X$ such that $U$ intersects all the nodes of $\eta$ disjoint from $X$, and every component of $G'-(V(U) \cup X)$ is adjacent in $G'$ to at most one vertex in $U$. 
By \cite[3.1]{rs} (taking $L=K_t$), there exist positive integers $\theta_0$ and $t_0$ (only depending on $t$ and hence only depending on $h$) such that for every graph $G'$ and every tangle in $G'$ of order at least $\theta_0$ not controlling a $K_t$-minor of $G'$, there exists a location $\LL$ in $G'$ contained in the tangle such that every torso of the design of $\LL$ is $K_{t_0}$-minor free.\footnote{More precisely, the location $\LL$ consists of the separations $(A,B)$ of $G'$, where $B=G-(V(S)-\overline{\Omega}))$ and $V(A)=V(S) \cup Z$, for each society $(S,\Omega)$ with $|\overline{\Omega}| \leq 3$ in the segregation of $G'-Z$ stated in \cite[3.1]{rs}, and consists of, for each society $(S,\Omega)$ with $|\overline{\Omega}|>3$ in the segregation of $G'-Z$ stated in \cite[3.1]{rs}, the separations $(A_1,B_1 \cup (G'-V(A_1))),...,(A_n,B_n \cup (G'-V(A_n)))$, where $(A_1,B_1),...,(A_n,B_n)$ are the separations of $S$ mentioned in \cite[3.2]{rs}. The property that the segregation is central implies that the location $\LL$ is contained in the tangle. As shown in the proof of \cite[1.3]{rs}, the torso of the design of $\LL$ can be made into an outgrowth by $\leq r_0$ $r_0$-rings of a graph that can be drawn in a surface in which $K_t$ cannot be drawn by deleting at most $r_0$ vertices, for some integer $r_0$ only depending on $t$ (and hence only depending on $h$). And it is well-known that such a graph is $K_{t_0}$-minor-free, for some integer $t_0$ only depending on $r_0$ (and hence only depending on $t$.}

Let $\D_1$ be the set of designs with $\D_1 = \{(Q,\C):$ every torso of $(Q,\C)$ is $K_{t_0}$-minor-free$\}$.
Let $\D_2$ be the set of designs with $\D_2 = \{(Q,\C):$ there exists $Z \subseteq V(Q)$ with $|Z| < 8h$ such that $Q-Z$ is bipartite and $|M-Z| \leq 1$ for every $M \in \C\}$.
Let $\D = \D_1 \cup \D_2$.

Let $\theta=\theta_0+8h$.
Now we prove that $\D$ is $\theta$-pervasive in $G$.
Let $G'$ be a subgraph of $G$.
Let $\T'$ be a tangle in $G'$ of order at least $\theta$.
If $\T'$ does not control a $K_t$-minor of $G'$, then by the second paragraph of the proof, there exists a location $\LL_1$ of $G'$ contained in $\T'$ such that every torso of the design of $\LL_1$ is $K_{t_0}$-minor free, so the design of $\LL_1$ is in $\D_1$.
So we may assume that $\T'$ controls a $K_t$-minor of $G'$, and we let $\alpha$ be the corresponding $K_t$-expansion.
Since $G'$ is a subgraph of $G$, $G'$ has no odd $K_h$-minor.
So by the second paragraph of the proof, there exist $X \subseteq V(G')$ with $|X| <8h$ and an induced bipartite subgraph $U$ of $G'-X$ such that $U$ intersects all nodes of $\alpha$ disjoint from $X$, and every component of $G'-(V(U) \cup X)$ is adjacent in $G'$ to at most one vertex in $U$. 
Hence there exists a location $\LL_2$ of $G'$ contained in $\T'$ such that the design of $\LL_2$ is $(G'[V(U) \cup X],\C)$, where $\C$ is a collection of subsets of $V(G')$ such that for every member $M$ of $\C$, $X \subseteq M$ and $|M-X| \leq 1$.
(We remark that $\LL_2$ is contained in $\T'$ since $t \geq 12h > |X|$.)
Note that $G'[V(U) \cup X]-X = U$ is bipartite.
So the design of $\LL_2$ is in $\D_2$.
Hence $\D$ is $\theta$-pervasive in $G$.

Let $\D_0^* = \{(Q,\C): |V(Q)| \leq 4\theta-3\}$ be a set of designs.
For each $i \in \{1,2\}$, let $\D^*_i$ be the set of designs such that $\D^*_i=\{(Q,\C):$ there exist $(Q',\C') \in \D_i$ and $Z \subseteq V(Q)$ with $|Z| \leq 3\theta-2$ such that $Q'=Q-Z$ and for every $M \in \C$ with $M \not \subseteq Z$, $M \cap V(Q') \in \C'\}$.
By \cite[2.1]{rs}, $G$ has a tree-decomposition $(B_x: x \in V(T))$ for some tree $T$ such that for every $x \in V(T)$, $(R_x,\C_x) \in \D_0^* \cup \D_1^* \cup \D_2^*$, where $R_x$ is the torso at $x$ and $\C_x$ is the frame at $x$.

Let $r=\max\{t_0+3\theta-2,4\theta-2\}$ and $\xi = 8h+3\theta-2$.
So $r$ and $\xi$ only depend on $h$.

Let $x \in V(T)$.
Let $R_x$ be the torso at $x$ and $\C_x$ be the frame at $x$.
Recall that $(R_x,\C_x) \in \D_0^* \cup \D_1^* \cup \D_2^*$.
If $(R_x,\C_x) \in \D^*_0 \cup \D^*_1$, then $R_x$ is $K_r$-minor free.
If $(R_x,\C_x) \in \D^*_2$, then there exists $Z \subseteq V(R_x)$ with $|Z| \leq \xi$ such that $R_x-Z$ is bipartite and $|M-Z| \leq 1$ for every $M \in \C$; note that in this case, since $|M-Z| \leq 1$ for every $M \in \C$, $R_x-Z$ is an induced subgraph of $G$.
This proves the theorem.
\end{pf}

\bigskip

Now we are ready to prove Theorem \ref{thm:odd_minor_intro}.
The following is a restatement.

\begin{theorem}
For every graph $H$, there exists a positive integer $c_H$ such that the following holds.
Let $S$ be a set of positive integers with $1 \in S$.
Let $G$ be a graph such that $H$ is not an odd minor of $G$.
	\begin{enumerate}
		\item If every induced bipartite subgraph of $G$ is $S$-achieved $k$-choosable, then $G$ is $S$-achieved $(k+c_H)$-choosable.
		\item If every induced bipartite subgraph of $G$ is properly $S$-achieved $k$-colorable, then $G$ is properly $S$-achieved $(k+c_H)$-colorable.
	\end{enumerate}
\end{theorem}

\begin{pf}
Let $H$ be a graph.
Let $r_H,\xi_H$ be the integers $r,\xi$ mentioned in Theorem \ref{thm:odd_structure}, respectively.
By \cite{m_minor}, there exists an integer $d_H$ such that every graph with no $K_{r_H}$-minor is $d_H$-degenerate.  
Define $c_H=2d_H+2\xi_H+2\max\{r_H-1,\xi_H+1\}$.

Let $S$ and $G$ be as stated in the theorem.
Let $\D_1$ be the set of designs with $\D_1=\{(Q,\C): Q$ is $K_{r_H}$-minor free and $\C$ is the set of all cliques of $Q\}$.
Let $\D_2$ be the set of designs with $\D_2=\{(Q,\C):$ there exists $Z \subseteq V(Q)$ with $|Z| \leq \xi_H$ such that $Q-Z$ is isomorphic to an induced bipartite subgraph of $G$, and for every $M \in \C$, $|M-Z| \leq 1\}$.
By Theorem \ref{thm:odd_structure}, $G$ has a tree-decomposition $\T$ over $\D_1 \cup \D_2$.
Note that the rank of $\D_1$ is at most $r_H-1$ and the rank of $\D_2$ is at most $\xi_H+1$.
So the adhesion of $\T$ is at most $\max\{r_H-1,\xi_H+1\}$.
For every $(Q,\C) \in \D_1$ and for every $\C$-compatible subgraph $H$ of $Q$, $Q-E(H)$ is $K_{r_H}$-minor free, so $Q-E(H)$ is $S$-achieved $(2d_H+1)$-choosable by Corollary \ref{cor:minor_cho}.
For every $(Q,\C) \in \D_2$ and for every $\C$-compatible subgraph $H$ of $Q$, there exists $Z \subseteq V(Q)$ with $|Z| \leq \xi_H$ such that $(Q-E(H))-Z$ is isomorphic to an induced subgraph of $G$, so for every set $X$ with $X \supseteq Z$, $(Q-E(H))-X$ is isomorphic to an induced subgraph of $G$, and hence by the assumption of the theorem, $(Q-E(H))-X$ is $S$-achieved $k$-choosable (and properly $S$-achieved $k$-colorable, respectively), and hence $Q-E(H)$ is $S$-achieved $(k+2\xi_H)$-choosable (and properly $S$-achieved $(k+2\xi_H)$-colorable, respectively) by Lemma \ref{lemma:deletion}.
Therefore, by Lemma \ref{lemma:design_torsos}, $G$ is $S$-achieved $(\max\{2d_H+1,k+2\xi_H\}+2\max\{r_H-1,\xi_H+1\})$-choosable (and properly $S$-achieved $(\max\{2d_H+1,k+2\xi_H\}+2\max\{r_H-1,\xi_H+1\})$-colorable, respectively).
Note that $k+c_H \geq \max\{2d_H+1,k+2\xi_H\}+2\max\{r_H-1,\xi_H+1\}$.
\end{pf}

\section{Vertex ordering} \label{sec:layer}

In this section we prove Theorem \ref{thm:ltw_intro}.
We need the following lemma.

\begin{lemma} \label{lemma:ordering}
Let $n$ be a positive integer.
Let $G$ be a graph on $n$ vertices.
For every $i \in [n]$, let $S_i$ be a subset of $V(G)$.
Let $v_1,v_2,...,v_n$ be an ordering of $V(G)$ such that for every $i \in [n]$, $N_G[v_i] \cap \{v_\ell: \ell \in [i]\} \subseteq S_i \subseteq \{v_\ell: \ell \in [i]\}$.
Let $w_1$ and $w_2$ be positive integers.
If $|S_i| \leq w_1$ and $|\bigcup_{j \in [n], v_i \in S_j}S_j \cap \{v_\ell: \ell \in [i]\}| \leq w_2$ for every $i \in [n]$, then $G$ is conflict-free $(w_1+w_2-1)$-choosable.
\end{lemma}

\begin{pf}
Let $L$ be a $(w_1+w_2-1)$-list-assignment of $G$.
For every $i \in [n]$, let $G_i=G[\{v_j: j \in [i]\}]$, and let $S_i^*=\bigcup_{j \in [n], v_i \in S_j}S_j$.
Note that $|S_i^* \cap V(G_i)| \leq w_2$ for every $i \in [n]$ by assumption.
We shall prove that for every $i \in [n]$, there exists a proper conflict-free $L|_{V(G_i)}$-coloring $\phi_i$ of $G_i$ such that for every $j \in [n]$, all vertices in $S_j \cap V(G_i)$ receive different colors in $\phi_i$.
Note that this statement implies the lemma since $G_n=G$.

We shall prove the claim by induction on $i$.
Since $w_1+w_2-1 \geq 1$, it obviously holds when $i=1$.
So we may assume that $i \geq 2$ and there exists a proper conflict-free $L|_{V(G_{i-1})}$-coloring $\phi_{i-1}$ of $G_{i-1}$ such that for every $j \in [n]$, all vertices in $S_j \cap V(G_{i-1})$ receive different colors in $\phi_{i-1}$. 

Since $\phi_{i-1}$ is a proper conflict-free $L|_{V(G_{i-1})}$-coloring of $G_{i-1}$, for every $j \in [i-1]$, either $N_{G_{i-1}}(v_j)=\emptyset$, or there exists $c_j$ such that $|\phi_{i-1}^{-1}(\{c_j\}) \cap N_{G_{i-1}}(v_j)|=1$.
Let $Z=\{\phi_{i-1}(v): v \in S^*_i \cap V(G_{i-1})\} \cup \{c_j: j \in [i-1], v_j \in S_i, N_{G_{i-1}}(v_j) \neq \emptyset\}$.
So $|Z| \leq |S^*_i \cap V(G_i)-\{v_i\}|+|S_i-\{v_i\}| \leq (w_2-1)+(w_1-1)$ since $v_i \in S_i \subseteq S^*_i$.
Since $L$ is a $(w_1+w_2-1)$-list-assignment of $G$, we can extend $\phi_{i-1}$ to a proper $L|_{V(G_i)}$-coloring $\phi_i$ of $G_i$ by further coloring $v_i$ by using an element in $L(v_i)-Z$.
We shall prove that $\phi_i$ is a desired coloring.

Suppose to the contrary that $\phi_i$ is not a conflict-free coloring of $G_i$.
So there exists $v \in V(G_i)$ such that $N_{G_i}(v) \neq \emptyset$ and no color in the image of $\phi_i$ appears exactly once in $N_{G_i}(v)$.
Since $\phi_{i-1}$ is conflict-free, $v_i \in N_{G_i}[v]$.
Hence $v \in N_{G_i}[v_i] \subseteq S_i$.
If $v \neq v_i$, then either $v_i$ is the unique neighbor of $v$ in $G_i$, or some color $c_v$ in $Z$ appears exactly once in $N_{G_{i-1}}(v)$.
For the former, $\phi_i(v_i)$ is a color appearing exactly once in $N_{G_i}(v)$, a contradiction; for the latter, since $\phi_i(v_i) \not\in Z$, $c_v$ is a color in the image of $\phi_i$ appearing exactly once in $N_{G_i}(v)$, a contradiction.
So $v=v_i$.
By assumption, $N_{G_i}(v) = N_G(v_i) \cap \{v_\ell: \ell \in [i-1]\} \subseteq S_i \cap V(G_{i-1})$.
So all vertices in $N_{G_i}(v)$ receive different colors in $\phi_{i-1}$ and hence in $\phi_i$.
Since $N_{G_i}(v) \neq \emptyset$, some color in the image of $\phi_i$ appears exactly once in $N_{G_i}(v)$, a contradiction.

So $\phi_i$ is a proper conflict-free $L|_{V(G_i)}$-coloring of $G_i$.
To prove this lemma, it suffices to show that for every $j \in [n]$, all vertices in $S_j \cap V(G_i)$ receive different colors in $\phi_i$.
Suppose to the contrary that there exists $j \in [n]$ such that at least two vertices in $S_j \cap V(G_i)$ receive the same color in $\phi_i$.
Since all vertices in $S_j \cap V(G_{i-1})$ receive different colors in $\phi_{i-1}$, we know that $v_i \in S_j$ and there exists $u \in S_j \cap V(G_i) \setminus \{v_i\} = S_j \cap V(G_{i-1})$ such that $\phi_i(v_i)=\phi_i(u)=\phi_{i-1}(u)$.
Since $v_i \in S_j$, $S_j \subseteq S_i^*$.
Since $u \in S_j \cap V(G_{i-1}) \subseteq S_i^* \cap V(G_{i-1})$, $\phi_i(v_i)=\phi_{i-1}(u) \in Z$, a contradiction.
This proves the lemma.
\end{pf}

\bigskip

A \emph{layering} of a graph $G$ is an ordered partition $(V_i: i \in {\mathbb N})$ of $V(G)$ such that, for each edge $vw$ of $G$, there exists an integer $i$ such that $\{v,w\}\subseteq V_i\cup V_{i+1}$. 
The \emph{layered treewidth} of a graph $G$ is $\min_{\V,\T}\max\{|V \cap B|: V \in \V, B \in \T\}$, where the minimum is over all layerings $\V$ of $G$ and all tree-decompositions $\T$ of $G$. 

Now we are ready to prove Theorem \ref{thm:ltw_intro}.
The following is a restatement.

\begin{theorem} \label{layer_cho}
Let $w$ be a positive integer.
If $G$ is a graph with layered treewidth at most $w$, then $G$ is conflict-free $(8w-1)$-choosable.
\end{theorem}

\begin{pf}
Since $G$ has layered treewidth at most $w$, there exist a layering $\V=(V_i: i \in {\mathbb N})$ of $G$ and a tree-decomposition $\T=(B_x: x \in V(T))$ of $G$ such that $|V \cap B| \leq w$ for every $V \in \V$ and $B \in \T$.

Let $r$ be a node of $T$.
We treat $T$ as a rooted tree rooted at $r$.
Let $x_1,x_2,...,x_{|V(T)|}$ be a breadth-first-search ordering of $V(T)$ with $x_1=r$.
For any two distinct nodes $a,b$ of $T$, define $a \prec_T b$ if $a=x_i$ and $b=x_j$ for some $i < j$.
Note that $\prec_T$ is a total order on $V(T)$.

For every $v \in V(G)$, let $x_v$ be the node $x$ of $T$ with $v \in B_{x}$ closest to $r$, and let $\ell_v$ be the integer such that $v \in V_{\ell_v}$.
Let $n=|V(G)|$.
Let $f$ be a bijection from $V(G)$ to $[n]$.
For any two vertices $u,v$ of $G$, define $u \preceq v$ if either 
	\begin{itemize}
		\item $x_u \prec_T x_v$, or 
		\item $x_u=x_v$ and $\ell_u < \ell_v$, or
		\item $x_u=x_v$ and $\ell_u=\ell_v$ and $f(u) \leq f(v)$.
	\end{itemize}
Then $\preceq$ is a total order on $V(G)$.

For every $i \in [n]$, 
	\begin{itemize}
		\item let $v_i$ be the $i$-th smallest vertex in $V(G)$ with respect to $\preceq$, 
		\item let $S_i=B_{x_{v_i}} \cap (\bigcup_{j=\ell_{v_i-1}}^{\ell_{v_i}+1} V_{j}) \cap \{v_j: j \in [i]\}$, and 
		\item let $Y_i=B_{x_{v_i}} \cap (\bigcup_{j=\ell_{v_i}-2}^{\ell_{v_i}+2} V_j) \cap \{v_j: j \in [i]\}$. 
	\end{itemize}
By the choice of $\V$ and $\T$, we know $|S_i| \leq 3w$ and $|Y_i| \leq 5w$.
Let $i \in [n]$.
Clearly, $S_i \subseteq \{v_\ell: \ell \in [i]\}$.

Suppose to the contrary that there exists $u \in N_G[v_i] \cap \{v_\ell: \ell \in [i]\} \setminus S_{i}$.
Since $\V$ is a layering of $G$, $u \in \bigcup_{j=\ell_{v_i}-1}^{\ell_{v_i}+1} V_j$. 
Since $u \in \{v_\ell: \ell \in [i]\}$, either $x_u \prec_T x_{v_i}$, or $x_u=x_{v_i}$ and $\ell_u \leq \ell_{v_i}$.
For the former, since $\T$ is a tree-decomposition, $u \in B_{x_{v_i}}$ by the definition of $x_{v_i}$, so $u \in S_i$, a contradiction; for the latter, $u \in B_{x_u}=B_{x_{v_i}}$, so $u \in S_i$, a contradiction.

Hence $N_G[v_i] \cap \{v_\ell: \ell \in [i]\} \subseteq S_{i} \subseteq \{v_\ell: \ell \in [i]\}$.

Suppose to the contrary that there exists $z \in \bigcup_{j \in [n], v_i \in S_j}S_j \cap \{v_\ell: \ell \in [i]\} \setminus Y_i$.
So there exists $j \in [n]$ such that $z \in S_{j}$ and $v_i \in S_{j}$.
Since $v_i \in S_j \subseteq B_{x_{v_j}}$, $x_{v_j}$ is a descendant of $x_{v_i}$ (including $x_{v_i}$) by the definition of $x_{v_i}$.
Since $v_i \in S_{j} \cap V_{\ell_{v_i}} \subseteq \bigcup_{k=\ell_{v_j}-1}^{\ell_{v_j}+1}V_k \cap V_{\ell_{v_i}}$, we have $|\ell_{v_j}-\ell_{v_i}| \leq 1$.
Since $z \in S_j \subseteq \bigcup_{k=\ell_{v_j}-1}^{\ell_{v_j}+1}V_k$, $|\ell_z-\ell_{v_j}| \leq 1$, so $|\ell_z-\ell_{v_i}| \leq 2$.
Since $z \in \{v_\ell: \ell \in [i]\}$, $x_z \preceq_T x_{v_i}$. 
Since $z \in S_j \subseteq B_{x_{v_j}}$ and $x_{v_j}$ is a descendant of $x_{v_i}$, we know that $x_z \preceq_T x_{v_i}$ implies $z \in B_{x_{v_i}}$.
Hence $z \in B_{x_{v_i}} \cap V_{\ell_z} \cap \{v_\ell: \ell \in [i]\} \subseteq B_{x_{v_i}} \cap \bigcup_{k=\ell_{v_i}-2}^{\ell_{v_i}+2}V_k \cap \{v_\ell: \ell \in [i]\} = Y_i$, a contradiction.

So $|\bigcup_{j \in [n], v_i \in S_j}S_j \cap \{v_\ell: \ell \in [i]\}| \leq |Y_i| \leq 5w$.
Therefore, $G$ is conflict-free $(3w+5w-1)$-choosable by Lemma \ref{lemma:ordering}.
\end{pf}

\bigskip

We remark that our proof for Theorem \ref{layer_cho} and its preparation Lemma \ref{lemma:ordering} was inspired by the proof of a result in \cite{dmo} about strong products of graphs.

One special case of graphs that have bounded layered treewidth is a strong product of a graph with bounded treewidth and a path.
For two graphs $A$ and $B$, the \emph{strong product} $A\boxtimes B$ of $A$ and $B$ is the graph with vertex-set $V(A)\times V(B)$ and that contains an edge with endpoints $(x_1,y_1)$ and $(x_2,y_2)$ if and only if
\begin{itemize} 
  \item $x_1x_2\in E(A)$ and $y_1=y_2$;
  \item $x_1=x_2$ and $y_1y_2\in E(B)$; or
  \item $x_1x_2\in E(A)$ and $y_1y_2\in E(B)$.
\end{itemize}
It is simple to show that if $A$ has treewidth at most $w$ and $B$ is a path, then $A \boxtimes B$ has layered treewidth at most $w+1$.
On the other hand, for every positive integer $w$, there exists a graph with layered treewidth 1 that cannot be written as a strong product of a graph with treewidth at most $w$ and a path \cite{bdjmw}.
Lemma \ref{lemma:ordering} can be used for graphs that are strong products of a bounded treewidth graph and a graph of bounded maximum degree by using a proof similar to the one of Theorem \ref{layer_cho}.
The following theorem implies an analogous result in \cite{dmo} stating that the same number of colors are enough for proper odd coloring of the same graph.

\begin{theorem}
Let $w$ and $d$ be nonnegative integers.
Let $H$ be a graph with treewidth at most $w$.
Let $Q$ be a graph with maximum degree at most $d$.
Then $H \boxtimes Q$ is conflict-free $((w+1)(d^2+d+2)-1)$-choosable.
\end{theorem}

\begin{pf}
Let $(B_x: x \in V(T))$ be a tree-decomposition of $H$ with width $w$.
Treat $T$ as a rooted tree with the root $r$.
Define $\prec_T$ to be a breadth-first-search ordering of $V(T)$ starting at $r$. 
For every $h \in V(H)$, let $x_h$ be the node $x$ of $T$ with $h \in B_x$ closest to $r$.
Let $\prec_Q$ be a total ordering of $V(Q)$.
Let $f$ be a bijection from $V(H)$ to $[n]$.
For any $(h_1,q_1),(h_2,q_2) \in V(G)$, we define $(h_1,q_1) \prec (h_2,q_2)$ if and only if either $x_{h_1} \prec_T x_{h_2}$, or $x_{h_1}=x_{h_2}$ and $q_1 \prec_Q q_2$, or $x_{h_1}=x_{h_2}$ and $q_1 = q_2$ and $f(h_1) < f(h_2)$.
Then $\prec$ is a total ordering of $V(G)$.

Let $n=|V(G)|$.
For every $i \in [n]$,
	\begin{itemize}
		\item let $v_i$ be the $i$-th smallest vertex of $V(G)$ with respect to $\prec$,
		\item let $h_i$ be the vertex of $H$ and $q_i$ be the vertex of $Q$ such that $v_i=(h_i,q_i)$,
		\item let $S_i = (B_{x_{h_i}} \times N_Q[q_i]) \cap \{v_\ell: \ell \in [i]\}$, and 
		\item let $Y_i = (B_{x_{h_i}} \times \{q \in V(Q):$ the distance in $Q$ between $q$ and $q_i$ is at most $2\}) \cap \{v_\ell: \ell \in [i]\}$.
	\end{itemize}
Then $|S_i| \leq (w+1)(d+1)$ and $|Y_i| \leq (w+1)(1+d+d(d-1))=(w+1)(d^2+1)$.
An argument similar as the proof of Theorem \ref{layer_cho} shows that $N_G[v_i] \cap \{v_\ell: \ell \in [i]\} \subseteq S_{i} \subseteq \{v_\ell: \ell \in [i]\}$ and $|\bigcup_{j \in [n], v_i \in S_j}S_j \cap \{v_\ell: \ell \in [i]\}| \leq |Y_i| \leq (w+1)(d^2+1)$.
Therefore, $G$ is conflict-free $((w+1)(d^2+d+2)-1)$-choosable by Lemma \ref{lemma:ordering}.
\end{pf}

\section{Concluding remarks} \label{sec:remarks}

When the writing of the first version of this paper \cite{l} was about to be completed, Hickingbotham \cite{h} announced a paper on arXiv. 
The main result in \cite{h} states that every graph $G$ is properly conflict-free $(2s_2-1)$-colorable, where $s_2$ is the 2-strong coloring number\footnote{The {\it $k$-strong coloring number} of a graph $G$ is the minimum $t$ such that there exists an ordering $v_1,v_2,...,v_{|V(G)|}$ of $V(G)$ such that for every $1 \leq i \leq |V(G)|$, there are at most $t$ indices $j \in [i-1]$ satisfying that there exists a path in $G$ from $v_i$ to $v_j$ with at most $k$ edges such that all internal vertices have indices greater than $i$.} of $G$. 
This result is essentially equivalent to (actually, slightly weaker than) Lemma \ref{lemma:ordering} of this paper (which is identical to its first version \cite[Lemma 5.1]{l}), with essentially the same proof, as taking $S_i$ in Lemma \ref{lemma:ordering} to be the union of $\{v_i\}$ and its left-neighborhood implies that every graph $G$ is conflict-free $(s_1+s_2-1)$-choosable, where $s_1$ is a number that is always at most $s_2$ and equals the 1-strong coloring number of $G$ for many natural graph classes, such as topological minor-closed classes. 
Lemma \ref{lemma:ordering} (and \cite[Lemma 5.1]{l}) was proved by the author of this paper before he knew \cite{h}. 
But he did not notice that it is equivalent to the formulation about strong coloring numbers until he saw \cite{h}.
Hickingbotham \cite{h} observed that combining the known result of Zhu \cite{z} about strong coloring numbers and his formulation of our Lemma \ref{lemma:ordering} immediately implies that every class with bounded expansion has bounded proper conflict-free chromatic number.

\begin{corollary}[\cite{h}] \label{cor:h_obs}
For every function $f: {\mathbb Z} \rightarrow {\mathbb Z}$, there exists an integer $c_f$ such that if $\F$ is a graph class such that for every graph $H$ in $\F$ and for every nonnegative integer $r$, every $r$-shallow minor\footnote{A graph $H'$ is an {\it $r$-shallow minor} of $H$ if $H'$ is isomorphic to a graph that can be obtained from a subgraph of $H$ by contracting disjoint subgraphs with radius at most $r$.} of $H$ has average degree at most $f(r)$, then $\chi_{\rm pcf}(G) \leq c_f$ for every graph $G$ in $\F$. 
\end{corollary}

Since our Lemma \ref{lemma:ordering} works for proper conflict-free list-coloring, by combing it with the same result of Zhu \cite{z}, we obtain that every graph in the class $\F$ stated in Corollary \ref{cor:h_obs} is actually conflict-free $c_f$-choosable.
Note that every graph with no subdivision of $H$ for any fixed graph $H$ lies in such a class $\F$.
So we obtain a proof of Theorem \ref{thm:topo_intro} different from the one stated in Section \ref{sec:topo}.
Note that the upper bounds for the number of colors in both proofs rely on implicit constants stated in results in the literature.
After a very quick and rough inspection for those implicit constants, the proof using clique-sums in Section \ref{sec:topo} seems giving a better upper bound, even though the proof that uses Lemma \ref{lemma:ordering} is conceptually simpler.
We keep both proofs because the proof in Section \ref{sec:topo} is a just simple application of our machinery for clique-sums, which leads to our result about odd minors (Theorem \ref{thm:odd_minor_intro}).
Note that Theorem \ref{thm:odd_minor_intro} does not follow from the list-coloring version of Corollary \ref{cor:h_obs} (and even Corollary \ref{cor:sparse} below), because odd minor-free graphs can be very dense and Theorem \ref{thm:odd_minor_intro} also involves induced subgraphs.

In fact, we can strengthen Corollary \ref{cor:h_obs} by combining Lemma \ref{lemma:ordering} with more results in the literature.

\begin{corollary} \label{cor:sparse}
For every positive integer $k$, there exists an integer $c$ such that if $G$ is a graph that is not conflict-free $c$-choosable, then $G$ contains a subgraph that is a 1-subdivision of a graph with average at least $k$ and contains a subgraph that is a 3-subdivision of a graph with chromatic number $k$.
\end{corollary}

\begin{pf}
Let $k$ be a positive integer.
By \cite[Lemma 4.5]{no}, for every integer $x$, there exists an integer $g(x) \geq x$ such that every graph with average degree at least $g(x)$ contains a $1$-subdivision of a graph with chromatic number $x$. 
Moreover, every $(\leq 1)$-subdivision of a graph with average degree at least $g(g(k)+1)+g(k)$ either contains a subgraph with average degree at least $g(g(k)+1)$ or contains a 1-subdivision of a graph with average degree at least $g(k)$; the former implies that it contains a subgraph that is a 1-subdivision of a graph that is minimal non-properly $g(k)$-colorable and hence has average degree at least $g(k)$; so either case implies that it contains a 1-subdivision of a graph with average degree at least $g(k) \geq k$ and hence contains a 3-subdivision of a graph with chromatic number $k$.
By \cite[p. 72]{no}, there exists an integer $c_1$ such that every graph with $\nabla_{1/2} \geq c_1$ contains an $(\leq 1)$-subdivision of a graph with average degree at least $g(g(k)+1)+g(k)$.
(We will not define $\nabla_{1/2}$ because we do not need the formal definition in this proof.)
By \cite[Corollary 3.5]{z}, there exists an integer $c_2$ such that every graph with 2-strong coloring number at least $c_2$ has $\nabla_{1/2} \geq c_1$.
By Lemma \ref{lemma:ordering}, there exists an integer $c$ such that every graph that is not conflict-free $c$-choosable has 2-strong coloring number at least $c_2$.
This proves the corollary.
\end{pf}

\bigskip

Corollary \ref{cor:sparse} strengthens Corollary \ref{cor:h_obs} not only to list-coloring but also to wider graph classes.
Corollary \ref{cor:sparse} is motivated by a potential coarse characterization of the graphs with large conflict-free choice number or $\chi_{{\rm pcf}}$ in terms of the subgraph relation. 
Recall that graphs with large chromatic number and 1-subdivision of graphs with large chromatic number are typical examples of graphs with large $\chi_{{\rm pcf}}$ (and hence with large conflict-free choice number).
Also, every graph with large chromatic number contains a subgraph with large average degree and hence contains a 1-subdivision of a graph with large chromatic number by \cite[Lemma 4.5]{no}.
So it is reasonable to ask whether 1-subdivision of graphs with large chromatic number are the only obstructions for having small conflict-free choice number or $\chi_{{\rm pcf}}$.

\begin{question} \label{que:opt_sub}
Does there exist a function $f$ such that for every integer $k$, every graph that is not conflict-free $f(k)$-choosable contains a subgraph that is a 1-subdivision of a graph with chromatic number at least $k$?
\end{question}

Note that Corollary \ref{cor:sparse} is evidence for a potential positive answer of Question \ref{que:opt_sub}.

Even though Question \ref{que:opt_sub} looks strong, it does not imply our result for odd minors (Theorem \ref{thm:odd_minor_intro}) since Theorem \ref{thm:odd_minor_intro} addresses induced subgraphs.
So it is reasonable to consider the following variant of Question \ref{que:opt_sub}.

\begin{question} \label{que:opt_induced}
Does there exist a function $f$ such that for every integer $k$, every graph that is not conflict-free $f(k)$-choosable (and not properly conflict-free $f(k)$-colorable, respectively) either has choice number at least $k$ (and chromatic number at least $k$, respectively) or contains an induced subgraph that is a 1-subdivision of a graph with choice number (and chromatic number, respectively) at least $k$?
\end{question}

A positive answer of Question \ref{que:opt_induced} would imply the case $S=\{1\}$ of Theorem \ref{thm:odd_minor_intro}.
To prove a positive answer of the coloring version of Question \ref{que:opt_induced}, it suffices to prove that every graph that is not properly conflict-free $f(k)$-colorable contains an induced subgraph that is an $(\leq 1)$-subdivision of a graph with chromatic number at least $k$ because of the following proposition.

\begin{proposition} \label{prop:at_most_to_exact}
Let $k$ be an integer.
Let $H$ be a graph with $\chi(H) \geq (k-1)^2+1$.
If $G$ is an $(\leq 1)$-subdivision of $H$, then there exists an induced subgraph $G'$ of $G$ such that either $\chi(G') \geq k$, or $G'$ is a 1-subdivision of a graph $H'$ with $\chi(H') \geq k$.
\end{proposition}

\begin{pf}
Let $H_1$ be the spanning subgraph of $H$ such that the edges of $H_1$ are exactly the edges of $H$ that are not subdivided in $G$.
Let $H_2 = H-E(H_1)$. 
Since $H_1$ is a subgraph of $G$, we may assume $\chi(H_1) \leq k-1$, for otherwise $\chi(G) \geq k$ and we are done.
Hence $V(H)$ can be partitioned into $k-1$ parts $V_1,V_2,...,V_{k-1}$ such that for every $i \in [k-1]$, some induced subgraph of $G$ is a 1-subdivision of $H[V_i]$.
Since $\chi(H) \leq \sum_{i=1}^{k-1} \chi(H[V_i])$, there exists $i^* \in [k-1]$ such that $\chi(H[V_{i^*}]) \geq \lceil \chi(H)/(k-1) \rceil \geq k$.
This proves the proposition.
\end{pf}

\bigskip

\noindent{\bf Acknowledgement:} The author thanks Vida Dujmovi\'{c}, Pat Morin, and Saeed Odak for some discussion about odd colorings and this paper.
The author also thanks anonymous reviewers for helpful comments, especially one of them pointed out that complete bipartite graphs provide a negative answer of an earlier version of Question \ref{que:opt_induced}.

\end{document}